\documentclass[12pt]{article}
\usepackage{amsthm,amsbsy,amsfonts,amssymb,amsmath,amscd}
\usepackage{latexsym}
\usepackage{fontenc}
\usepackage[mathscr]{euscript}

\usepackage[dvips]{epsfig}
\DeclareMathAlphabet{\scrb}{U}{eus}{b}{n}
\theoremstyle{plain}

\newtheorem{Thm}{Theorem}[section]

\newtheorem{Prop}{Proposition}[section]
\newtheorem{lem}{Lemma}[section]
\newtheorem{cor}{Corollary}[section]
\theoremstyle{remark}

\newtheorem{Thme}[Prop]{\bf Theorem}

\begin{document}
\def\bcw{\mathbin{\bigcirc\mkern-15mu\wedge}}

\title{{\bf Conformal Geometry on Four Manifolds}}  

\author{{\bf ICM Emmy Noether Lecture} \\  {\bf ICM 2018, Rio de Janeiro, Brazil}  \\  \\ {\bf Sun-Yung Alice Chang }\footnote{Research of Chang is supported
in part by NSF Grant DMS-1509505}}
\date{}
\maketitle

\setcounter{section}{-1}
\setcounter{equation}{00}
\noindent 
{\bf \S 0. Introduction}
\vskip.1in
This is the lecture notes for the author's Emmy Noether lecture at 2018, ICM, Rio de Janeiro, Brazil. It is a great honor for
the author to be invited to give the lecture. 

In the lecture notes, the author will survey the development of conformal geometry on four dimensional manifolds.
The topic she chooses is one on which she has been involved in the past twenty or more years:  the study of the integral conformal invariants on 4-manifolds and geometric applications.  
The development was heavily influenced by many earlier pioneer works; recent progress in conformal geometry has also been made
in many different directions, here we will only present some slices of the development. 

The notes is organized as follows. 

In section 1, we briefly describe the prescribing Gaussian curvature problem on compact surfaces and the Yamabe problem
on $n$-manifolds for $n \geq 3$;  in both cases some second order PDE have played important roles. 

In section 2, we introduce the quadratic 
curvature polynomial $\sigma_2$ on compact closed 4-manifolds, which appears as part of the integrand of the Gauss-Bonnet-Chern formula. We discuss its algebraic structure, its connection to the 4-th order Paneitz operator $P_4$ and its associated 4-th order $Q$ curvature. We also discuss some variational approach to study the curvature and as a geometric application, results  
to characterize the diffeomorphism type of $(S^4, g_c)$ and $({\mathbb {CP}}^2, g_{FS})$ in terms of the size of the conformally
invariant quantity: the integral of $\sigma_2$ over the manifold.

In section 3, we extend our discussion to  compact
4-manifolds with boundary and introduce a third order pseudo-differential operator $P_3$ and 3-order curvature $T$ on 
the boundary of the manifolds. 

In section 4, we shift our attention to the class of conformally compact Einstein (abbreviated as CCE)  four-manifolds.  We survey some recent 
research on the problem of ``filling in" a given 3-dimensional manifold as the conformal infinity of a CCE manifold.
We relate the concept of "renormalized" volume in this setting again to the integral of $\sigma_2$.

In section 5, we discuss some partial results on a compactness problem on CCE manifolds. We believe the compactness results are the 
key steps toward an existence theory for CCE manifolds.

The author is fortunate to have many long-term close collaborators, who have greatly contributed to the development of the research described in this article -- some more than the author. Among them Matthew Gursky, Jie Qing, Paul Yang and more
recently Yuxin Ge. She would like to take the chance to express her deep gratitude toward them, for the fruitful collaborations and for the 
friendships.

\vskip .4in

\setcounter{section}{1}
\noindent
{\bf \S 1. \,\, Prescribing Gaussian curvature on compact
surfaces and the Yamabe \,\,\,\,\, problem}
\setcounter{equation}{00}
\vskip .1in
In this section we will describe some second order elliptic equations
which have played important roles in conformal geometry.

On a compact surface $(M, g)$ with a Riemannian metric $g$, a natural
curvature invariant associated with the Laplace operator $\Delta = \Delta
_g$ is the Gaussian curvature $K=K_g$. Under the conformal change of
metric $g_w = e^{2w} g$, we have
\begin{equation}
- \Delta w \, + \, K = K_{g_w} e^{2w} \,\, on \,\, M.
\end{equation}
The classical
uniformization theorem to classify compact closed surfaces can be viewed
as finding solution of equation (1.1) with $K_{g_w} \equiv -1$, $0$, or
$1$
according to the sign of $\int K_g dv_g$. Recall that the
Gauss-Bonnet theorem states
\begin{equation}
2 \pi \, \chi(M) \, = \,  \int_M  K_{g_w} \, dv_{g_w}, 
\end{equation}
where $\chi(M)$ is the Euler characteristic of $M$, a topological
invariant. 
The variational functional  with
(1.1) as Euler equation for $K_{g_w} = \text{constant} $  is thus given by Moser's functional (\cite{M71}, \cite{M71b})
\begin{equation}
J_g [w] = \int_M  |\nabla w |^2 dv_g + 2 \int_M  K w dv_g -\left(\int_M K dv_g\right) \log
\frac { \int_M dv_{g_w} }{\int_M dv_g}.
\end{equation}

There is another geometric meaning of the functional $J$ which influences the later development of the field,
that is the formula of Polyakov \cite{Po} 
%%(see also Okikiolu \cite {Ok-1})
\begin{equation}
\label{logdet}
J_g [w]=  12 \pi \,\,\, \log \,\,\left(\frac{\det \,\, {(-\Delta)} _g}{\det \,\, {(-\Delta)} _{g_w}}\right)
\end{equation}
for metrics $g_w$ with the same volume as $g$;
where the determinant of the Laplacian $\det \,\, \Delta _g$ is defined by Ray-Singer via
the ``regularized'' zeta function. 

In \cite{On}, (see also Hong \cite{Ho}), Onofri established the sharp
inequality that on the 2-sphere
$J [w] \geq 0$  and $ J [w] = 0 $ precisely for conformal factors $w$ of the form
$e^{2w}g_0= T^*g_0$ where $T$ is a Mobius transformation of the
2-sphere. Later Osgood-Phillips-Sarnak (\cite{OPS88}, \cite{OPS88b}) arrived at the same sharp inequality in their
study
of heights of the Laplacian. This inequality also
plays an important role in their proof of the $C^{\infty}$
compactness of isospectral metrics on compact surfaces. 

%The formula of Polyakov-Ray-Singer has been generalized to manifolds of dimension
%greater than two in many different
%settings; one of which we will discuss in section 2 below.
%There is also a general study of extremal metrics for $ det \,\, \Delta_g
%$ or $ det \,\, L_g$
%for metrics $g$  in the same conformal class
%with a fixed volume
%or for all metrics with a fixed volume. (\cite {Be}, \cite {BCY}, \cite {Br-3}, \cite {R},
%\cite {Ok-2}). A special case of the remarkable results of
%Okikiolu (\cite {Ok-2})
%is that among all metrics with the same volume as the standard metric on
% the 3-sphere, the standard canonical metric is a local maximum for the
%functional $det \,\,\Delta_g$.

On manifolds $(M^n, g)$ for n greater than two,  the conformal
Laplacian $L_g$ is defined as
$L_g = - \Delta_g + c_n  R_g$ where $c_n = \frac {n-2} {4 (n-1)}$, and $R_g$ denotes the scalar
curvature of the metric $g$. An analogue of equation (1.1) is the equation, commonly referred to as the Yamabe equation (\ref{scalar}), which relates
the scalar
curvature under conformal change of metric to the background metric. In this case, it is
convenient to denote the conformal metric as $ \hat g = u ^{ \frac {4}{n-2}} g$ for some positive function $u$, then the equation becomes
\begin{equation}
\label{scalar}
L_g u \,\, = c_n \,\, \hat R \,\, u^{\frac {n+2}{n-2}}.
\end{equation}
The famous Yamabe problem to solve (\ref{scalar}) with $ \hat R $ a constant has been
 settled by Yamabe \cite{Y60}, Trudinger \cite{Tr68}, Aubin \cite{Au76} and
Schoen \cite{S84}. The corresponding problem to prescribe scalar
curvature has been intensively studied in the past decades by different
groups of
mathematicians, we will not be able to survey all the results here. We will only point out that 
in this case the study of $\hat R = c $, where $c$ is a constant, over class of metrics $\hat g$ in the conformal class $[g]$ 
with the same volume as $g$, is a variational problem with respect to the functional $F_g [u] = \int_{M^n}  R_{\hat g} dv_{\hat g}$
for any $n \geq 3$.  Again the sign of the constant $c$ agrees with the sign of the Yamabe invariant
\begin{equation}
\label{yam}
Y (M, g) : =  \inf _{\hat g \in  [g]} \frac {\int_M  R_{\hat g} dv_{\hat g}} {{(vol \,{\hat g} ) }^ { \frac {n-2}{n}}}. 
\end{equation}

\def\intl{\int\kern-9pt\hbox{$\backslach$}}

\vskip .2in
\noindent
{\bf \S 2. $\sigma_2$ curvature on 4-manifold}
\setcounter{equation}{00}

\setcounter{section}{2}
\setcounter{equation}{00}
\setcounter{Prop}{00}
\setcounter{Thm}{00}
\vskip .2in

%Define curvature, Ricci, 
\noindent {\bf \S 2a definition and structure of $\sigma_2$.}
\vskip .1in

We now introduce an integral conformal invariant 
which plays a crucial role in this paper,  namely the integral of $\sigma_2$  curvature on 
four-manifolds.
\vskip .1in
To do so, we first recall the Gauss-Bonnet-Chern formula on closed compact manifold $(M, g)$
 of dimension four:

\begin{equation}
8\pi^2 \chi (M)= \int _M \frac {1}{4} |W_g|^2 dv_g+ \int_M \frac {1}{6} (R_g^2 - 3 |Ric_g|^2)  dv_g,
\end{equation}

\noindent where $ \chi (M) $ denotes the Euler characteristic of $M$,  $W_g$ denotes the Weyl curvature, $R_g$ the scalar curvature and 
$Ric_g$ the Ricci curvature of the metric g.

%% May need to define Weyl etc)
\vskip .1in
In general, the Weyl curvature measures the obstruction to being conformally flat.  More precisely, for a manifold of dimension greater or equals to four,
$W_g $ vanishes in a neighborhood of a point if an only if the metric is locally conformal to a Euclidean metric; i.e., there are local coordinates such that $g = e^{2w} |dx|^2 $
for some function $w$.  Thus for example, the standard round metric $g_c$ on the sphere $S^n$ has $W_{g_c} \equiv 0$.
%In contrast, the Ricci curvature and the scalar curvature 
%measures the curvedness of the
%manifold, on $(S^n, g_c)$, $Ric_{g_c}  \equiv (n-1) $, $R_{g_c} \equiv n(n-1)$. 

\vskip .1in
In terms of conformal geometry, what is relevant to us is that Weyl curvature is a {\em pointwise}  conformal invariant,
in the sense that
under conformal change of metric $g_w = e^{2w} g$, $|W_{g_w}| = e^{-2w} |W_g| $, thus on 4-manifold 
 $|W_{g_w} |^2 dv_{g_w} = |W_g|^2 dv_{g} $;
this implies in particular that the first term in the Gauss-Bonnet-Chern formula above 
$$ g \rightarrow \int_M  |W|_g^2 dv_g $$
is conformally invariant.

\vskip .1in 
For reason which will be justified later below, we denote 
\begin{equation}
 \sigma_2 (g) = \frac {1}{6} (R_g^2 - 3 |Ric_g|^2) 
 \end{equation}
and draw the conclusion from the above discussion of the Gauss-Bonnet-Chern formula that
$$ g \rightarrow \int_M \sigma_2 (g) dv_g $$
is also an integral conformal invariant. This is the fundamental conformal invariant which will be studied in this lecture notes. 
We begin by justifying the name of ``$\sigma_2$" curvature.  
\vskip .1in
On manifolds of dimensions greater than two, the Riemannian curvature tensor $Rm$ can be decomposed 
into the different components.  From the perspective of conformal geometry, a natural basis is the 
Weyl tensor $W$, and the Schouten tensor, defined by
$$A_g =Ric_g- \frac{R}{2(n-1)}g.$$
The curvature tensor can be decomposed as
$$
{Rm}_g =W_g  \oplus  \frac{1}{n-2} A _g\bcw  g.
$$
Under conformal change of metrics $g_{w}=e^{2w}g$, since the Weyl tensor $W$ transforms by scaling,
only the Schouten tensor depends on the derivatives of the
conformal factor. It is thus natural to consider $\sigma_k(A_g)$, the k-th symmetric function of the eigenvalues
of the Schouten tensor $A_g$, as curvature invariants of
the conformal metrics.
 \vskip .1in 
When $k=1$, $\sigma_1(A_g) = Tr_g \,\,  A_g =  \frac{n-2}{ 2(n-1)} R_{g}$, 
so the $\sigma_1$-curvature is a dimensional multiple of the scalar curvature. 

When $ k =2 $, $\sigma_2 (A_g) = \sum_{i < j}  \lambda_{i} \lambda_{j} =  \frac {1}{2} (|Tr_g \,\, A_g|^2 - |A_g|^2)$,
where the $ \lambda $s are the eigenvalues of the tensor $A_g$. For a manifold of dimension 4,
we have 
$$ \sigma_2 (g) = \sigma_2 (A_g)  = \frac {1}{6} (R_g^2 - 3 |Ric_g|^2).$$
 
When $k =n$, $\sigma_n (A_g) = determinant \,\, of \,\, A_g$, an equation
of Monge-Amp\`ere type.
\vskip .1in 

In view of the Yamabe problem, it is natural to ask the question under what condition can one find a 
metric $g_w$ in the conformal class of $g$, which solves the equation 
\begin{equation}
\label{sig2}
{\sigma_2} (A_{g_w} ) = \, constant .\, 
\end{equation}
To do so, we
first observe that as a differential invariant of the conformal factor $w$, $\sigma_k(A_{g_w})$ is a
fully nonlinear expression involving the Hessian and the gradient of
the conformal factor $w$. We have
%\begin{equation}
$$
A_{g_w} \, = \,  (n-2) \{- \nabla ^2w +  dw\otimes dw- \frac{|\nabla w|^2}{2} \}
+ A_g.
$$
%\end{equation}
%The equation
%\begin{equation}
%\label{sig}
%\sigma _k(A_{g_w})= 1
%\end{equation}
%is a fully nonlinear version of the Yamabe equation.
% For example,
% when $k =1$, $\sigma_1(A_g) = \frac{n-2}{ 2(n-1)} R_{g}$, where
%equation \ref{sig}
%is the Yamabe equation which we have discussed in section 1. 
% When
% $ k =2 $, $\sigma_2 (A_g) = \frac {1}{2} (|Trace \,\, A_g|^2 - |A_g|^2 = \frac{n}{8 (n-1)} R^2 - \frac{1}{2} |Ric|^2 $. 
%In the case when $k =n$,
% $\sigma_n (A_g) = determinant \,\, of \,\, A_g$, an equation
%of Monge-Ampere type. 
\vskip .1in 

To illustrate that (\ref{sig2})  is a fully non-linear equation, we have when $n=4$,
\begin{equation}
\label{sig} 
 \begin{aligned}
\sigma_2(A_{g_w}) e^{4w} \, =&\, \sigma_2 (A_g) \, + 2 ((\Delta w)^2 
- \, |\nabla ^2 w|^2 \\ 
+& \,  \Delta w |\nabla w |^2  + (\nabla w,\nabla |\nabla w|^2) ) 
\,\, \\
+& \text {lower order terms} . 
\end{aligned}
\end{equation}
where all derivative are taken with respect to the $g$ metric.

For a symmetric $n\times n$ matrix $M$, we say $M \in \Gamma _k^+$
in the sense of G{\aa}rding  (\cite {Ga59})
if
$\sigma _k(M) >0$ and $M$ may be joined to the identity matrix by a path
consisting entirely of matrices $M_t$ such that $\sigma _k(M_t) >0$.
There is a rich literature concerning the equation
\begin{equation}
\label{si}
\sigma _k(\nabla ^2 u)\, = \, f ,
\end{equation}
for a positive function $f$, which is beyond the scope of this article to cover. 
Here we willl only note that when $k=2$,
\begin{equation}
\label{si2}
\sigma _2(\nabla ^2 u)\, =\,  (\Delta u)^2 - |\nabla^2 u|^2  \, =  \, f ,
\end{equation}
for a positive function $f$.  We remark the leading term of equations (\ref{si2}) and (\ref{sig}) agree. 
\vskip .1in

%In the case when $M = ( \nabla^2 u )$ for convex
%functions $u$ defined on the Euclidean domains, regularity theory
%for equations of $\sigma_k(M)$ has been well established for
%$M \in \Gamma _k^+$ for Dirichlet boundary value problems by
%Caffarelli-Nirenberg-Spruck (\cite{CNS85}); for a more general
%class of fully non-linear elliptic equations not
%necessarily of divergence form by Krylov (\cite{Kr82}),
%Evans (\cite{Ev82}) and for Monge-Ampere equations by Pogorelov (\cite{Pog71})
%and by Caffarelli (\cite{Ca90}). 
%The Monge-Ampere equation for prescribing the
%Gauss-Kronecker curvature for convex hypersurfaces has been studied by 
%Guan-Spruck \cite{GS}. 
%Some of the techniques in these work
%can be modified to study equation  (\ref{sig}) on manifolds. However there are
%features of the equation (\ref{sig}) that are distinct from the equation  (\ref{si}) when $f\equiv 1$.
\vskip .1in

We now discuss a variational approach to study the equation (\ref{sig2})  for $\sigma_2$ curvature .
\vskip .1in
Recall in section 1 we have mentioned that the functional 
$  {\cal F}(g) \,:= \, \int_{M^n}  R_g dv_g $ is variational in the sense that when $n \geq 3$ and when one varies $g$ in the same conformal
class of metrics with fixed volume, the critical metric when attained satisfies $R_{\hat g} \equiv \text{constant}$;  while when $n=2$, this is no longer true with $R_g$ replaced by $K_g$ and one needs to replace
the functional ${\cal F}(g) $ by the Moser's functional $J_g$. 
\vskip .1in
Parallel phenomenon happens when one studies the $\sigma_2 $ curvature. It turns out that when $n >2$ and $n \ne 4$, 
the functional  $ {{\cal F}_2} (g)  \, := \, \int_{M^n}  \sigma_2 (g) dv_g $ is variational when one varies $g$ in the same conformal
class of metrics with fixed volume, while this is no longer true when $n =4$. In section 2a below, we will describe
a variational approach to study equation (\ref{sig2}) in dimension 4 and the corresponding Moser's functional. Before we do so,
we would like to end the discussion of this section by quoting a result of Gursky-Viaclovsky \cite{GV01}. 
\vskip .1in

In dimension 3, one can capture all metrics with constant
sectional curvature (i.e. space forms) through the study of $\sigma_2$.

\begin{Thme} (\cite{GV01})  On a compact 3-manifold, for any Riemannian metric $g$,\, denote
${{\cal F}_2}(g)  = \int_M \sigma_2 (A_g) dv_g$.  Then a metric $g$
with $ {{\cal F}_2 }(g ) \geq 0$ is critical
for the functional ${\cal F}_2 $ restricted to class of metrics with volume
one if and only if $g$ has constant sectional curvature.
\end{Thme}

\vskip .2in

\noindent {\bf \S 2b. 4-th order Paneitz operator, Q-curvature}
\vskip .1in
We now describe the rather surprising link between a 4-th order linear operator $P_4$, 
its associated curvature invariant $Q_4$,  and the $\sigma_2$-curvature.

We first recall on $(M^n, g)$, $ n \geq 3$, the second order conformal Laplacian operator $L= -\Delta + \frac{n-2} {4(n-1)}  R$
transforms under the conformal change of metric $\hat g=u^{\frac{4}{n-2}}g$, as
\begin{equation}
L_{\hat g}(\varphi) = u^{-\frac{n+2}{n-2}}L_g(u \, \varphi \, ) \,\, \text{\rm for \,\, all}\quad \varphi \in C^\infty (M^4). \,
\end{equation}

% In general, we call a metrically defined operator $A$ conformally covariant
% of bidegree $(a, b)$,  if under the conformal change of metric
% $g_w = e^{2\omega }g$, the pair of corresponding operators
% $A_\omega $ and $A$ are related by
%$$
%A_\omega (\varphi ) = e^{-b\omega } A(e^{a\omega }\varphi )\quad
%\text{\rm for all}\quad \varphi \in C^\infty (M^n) \, .
%$$
%Note that in this notation, the conformal Laplacian operator is conformal
%of bidegree $(\frac{n-2}{2}, \frac{n+2}{2})$.

There are many operators besides the Laplacian $\Delta $
 on compact surfaces and the conformal Laplacian $L$ on general compact
 manifold of dimension greater than two  which have the conformal covariance
 property. One class of such operators of order 4 was studied by Paneitz (\cite{Pa83}, see also \cite{ES85})
 defined on $(M^n, g)$ when $n > 2$; which we call the conformal Paneitz operator:
\begin{equation}
P_4^n= \Delta ^2 + \delta \left( a_n Rg + b_n \text{\rm Ric} \right)
 d  + \frac{n-4}{2} Q_4^n
 \end{equation}
 and
 \begin{equation}
 Q_4^n = c_n |Ric |^2 + d_n R^2 - \frac{1}{2(n-1)}\Delta R,
 \end{equation}
 where $a_n$, $b_n$, $c_n$ and $d_n$ are some dimensional constants, 
 $\delta $ denotes the divergence, $d$ the deRham differential.  
%$$a_n = \frac {(n-2)^2 + 4}{2(n-1)(n-2)},\;
% b_n= - \frac 4{n-2},\;
 %c_n = -\frac{2}{(n-2)^2},\;
 %d_n =\frac{n^3-4n^2+16n-16}{8(n-1)^2(n-2)^2}.
%$$

The conformal Paneitz operator is conformally covariant. In this case, we write the conformal metric as $ \hat g=u^{\frac{4}{n-4}}g$ for some positive function $u$, then for $ n \ne 4$, 
 \begin{equation}
 {(P_4^n)}_{\hat g}(\varphi) \,  = \,  u^{-\frac{n+4}{n-4}} {(P_4^n)}_g (u \,\varphi\,) \,\, \, \text{ for  all} \,\,  \varphi \in C^\infty (M^4). \, 
 \end{equation}

Properties of $(P_4^n, Q_4^n)$ have been intensively studied in recent years, with many surprisingly strong results. We refer
the readers to the recent articles (\cite{CaY13}, \cite{GuM15}, \cite{GuHL16}, \cite{HY15}, \cite{HY16} and beyond).  \\

Notice that when $n$ is not equal to 4, we have $P_4^n (1) = \frac{n-4}{2} Q_4^n$, while when $n=4$, one
does not read $Q_4^4$ from $P_4^4$; it was pointed out by T. Branson that nevertheless 
both $P :=  P_4^4 $ and
$Q := Q_4^4$ are well defined (which we named as {\it {Branson's Q-curvature}}):

 $$
 P\varphi :=  \Delta ^2 \varphi + \delta
 \left( \frac 23 Rg - 2 \text{\rm Ric}\right) d \varphi.
 $$
 The {\it Paneitz} operator $P$
 is conformally covariant of bidegree
 $(0, 4)$ on 4-manifolds, i.e.
 $$
 P_{g_w} (\varphi ) = e^{-4\omega } P_g (\varphi )\quad
 \text{\rm for all}\quad \varphi \in C^\infty (M^4) \, .
 $$
The $Q$ curvature associated with $P$ is defined as
\begin{equation}
\label{Q4} 
2 Q_g \, = \,  - \frac{1}{6} \Delta R_g  \, + \, \frac {1}{6} (R_g^2 - 3 |Ric_g|^2)).
\end{equation} 
Thus the relation between $Q$ curvature and $\sigma_2$ curvature is 
\begin{equation}
\label{Q4s2}
2 Q_g \, = \, - \frac{1}{6} \Delta R_g + \sigma_2 (A_g).
\end{equation}
The relation between $P$ and $Q$ curvature on manifolds of dimension four is like that of the Laplace operator $- \Delta$ and the Gaussian curvature $K$
on compact surfaces.   
\begin{equation}
\label{PQ}
 P_g w + 2 Q_g \, = \, 2 Q_{g_w} e^{4w}.
 \end{equation}

Due to the pointwise relationship between $Q$ and $\sigma_2$  in (\ref{Q4s2}), notice on compact closed 4-manifolds $M$, 
$ \int_M  \Delta R_g  dv_g = 0$, thus
$ \int_M Q_g dv_g = \frac{1}{2} \int \sigma_2 (A_g) dv_g $ is also an integral conformal invariant. But $Q$ curvature has the advantage
of being prescribed by the linear operator $P$ in the equation (\ref{PQ}), which is easier to study. 
We remark the situation is different
on compact 4-manifolds with boundary which we will discuss in later part of this article \\

Following Moser, the functional to study constant $Q_{g_w}$ metric with $g_w \in [g]$ is given by 

$$ 
II[w] =\langle
Pw,w\rangle+4\int\,Qwdv-\left(\int\,Qdv\right)\,\log\,
\frac {\int  e^{4w}dv}{ \int dv} .
$$
In view of the relation \label{Q4} between $Q$ and $\sigma_2$ curvatures, if one consider the variational functional III
whose Euler equation is $ \Delta R = \, constant \, $, 
 
$$
III[w] = \frac {1}{3} \left ( \int R_{g_w}^2 dv_{g_w} - \int R^2 dv \right ),
$$

\noindent and define 

$$
{\cal F}  [w] = II [w] - \frac {1}{12} III [w], 
$$
we draw the conclusion:
\begin{Prop}  (\cite{CY03}, \cite{CGY02}, see also \cite{BV04} for an alternative approach.)\\
${\cal F} $ is the Lagrangian functional for $\sigma_2$ curvature. 
\end{Prop}  
\vskip .1in
We remark that the search for the functional $\cal F$ above was originally motivated by the study of some other variational formulas,  e.g.  the variation of quotients of log determinant of conformal
Laplacian operators under conformal change of metrics on 4-manifolds,  analogues to that  of the Polyakov formula (\ref{logdet}) on compact surfaces .  We refer the readers
to articles (\cite{BO91}, \cite{BCY92}, \cite{CGY02}, \cite{Ok01}) on this topic.

We also remark that on compact manifold of dimension $n$, there is a general class of conformal covariant operators $P_{2k}$ of order $2k$, 
for all integers $k$ with $ 2 k \leq n$. This is the well-known class of GJMS operators \cite{GJMS92}, where $P_2$ coincides with the conformal 
Laplace operator $L$ and where $P_4$
coincides with the 4-th order Paneitz operator. GJMS operators have played important roles in many recent developments in conformal geometry.
 
\vskip .2in

\noindent {\bf \S 2c  Some properties of the  conformal invariant $\int \sigma_2$}

\vskip .1in

We now concentrate on the class of compact, closed four manifolds which allow a Riemannian metric $g$ in the class $\cal A$, where    

$$ {\cal A} :=  \left \{g |, Y(M, g) >0, \,\, \int_{M} \sigma_2 (A_g) dv_g >0  \right\}.
$$

Notice that for closed four manifolds, it follows from equation (\ref{Q4s2}) that 
$$  2 \int_M Q_g dv_g \, = \, \int_{M} \sigma_2 (A_g) dv_g, $$
so in the definition of $\cal A$, we could also use $Q$ instead of $\sigma_2$. 
\vskip .1in
We recall some important properties of metrics in ${\it A}$.

\begin{Thme} (Chang-Gursky-Yang \cite{CY95}, \cite{Gu99}, \cite{CGY02}) \newline

1. If $Y(M,g)>0$, then $\int_M \sigma_2 (g) dv_g \leq 16 \pi^2$, equality holds if and only if $(M, g)$ is 
conformally equivalent to $(S^4, g_c)$.
\vskip .1in
2. $g \in {\cal A  }$, then $P \geq 0$ with $ \text{kernel} \, (P)$ consists of constants; it follows there exists some $g _w \in [g]$ with
$Q_{g_w} = constant $ and $R_{g_w} >0$.
\vskip .1in
3. $g \in {\cal A }$, then there exists some $g_w \in [g]$ with $ \sigma_2 (A_{g_w} ) >0 $ and $R_{g_w} >0$; i. e. $g_w $ exists in the positive two
cone  $\Gamma_2^{+}$ of $A_g$ in the sense of G{\aa}rding \cite{Ga59}. We remark  $ g  \in \Gamma_2^{+}$ 
implies $Ric_g >0$, as a consequence the first betti number $b_1$ of $M$ is zero.
\vskip .1in
4.  $g \in {\cal A}$, then there exists some $g_w \in [g]$ with $ \sigma_2 (A_{g_w} ) =1 $ and $R_{g_w} >0. $
\vskip .1in
5. When $(M, g)$ is not conformally equivalent to $(S^4, g_c) $ and  $g \in {\cal A}$, then for any positive smooth function $f$ defined on $M$,  there exists some $g_w \in [g]$ with $ \sigma_2 (A_{g_w} ) =f $ and  $R_{g_w} >0. $
\vskip .1in
\end{Thme}
\vskip .1in
We remark that techniques for solving the $ \sigma_2 $ curvature  equation can be modified to solve 
the equation $ \sigma_2 =1 + c |W|^2 $ for some constant $c$, which is the equation we will use later in the proof of 
the theorems in section 2d. \\

As a consequence of above theorem we have \\
\begin{cor}
On $(M^4, g)$, $g \in {\cal A}$ if and only if there exists some $g_w \in [g]$
with $g_w \in \Gamma_2^{+}$.
\end{cor}
\vskip .1in
A significant result in recent years is the following ``uniqueness" result.\\

\begin{Thme} (Gursky-Steets \cite{GuSt16}) \\
Suppose $(M^4,g)$ is not conformal to $(S^4, g_c)$ and $g \in \Gamma_2^{+}$, then $g_w \in [g] $ with $g_w \in \Gamma_2^{+}$ and with $ \sigma_2 (A_{g_w} ) =1 $ is unique.
\end{Thme}

The result was established by constructing some norm for metrics in $\Gamma_2^{+}$, with respect to which the functional
$\cal F$ is convex. The result is surprising in contrast with the famous example of R. Schoen {\cite{S87}} where he showed that on $(S^1 \times S^n,  g_{prod}) $, where $n \geq 2$,  the class of constant scalar curvature metrics (with the same volume) is not unique.

%\end{Thme}

\vskip .2in
\noindent {\bf  \S 2d   Diffeomorphism type}
\vskip .1in

In terms of geometric application, this circle of ideas may be applied to
characterize the diffeomorphism type of manifolds in terms of the
the relative size of the conformal invariant $\int \sigma _2(A_g)dV_g$
compared with the Euler number of the underlying manifold,  or equivalently the relative size of the two
integral conformal invariants $\int \sigma _2 (A_g) dv_g$ and $\int ||W||_g^2 dv_g$.  

Note: In the following, we will view the Weyl tensor as an endomorphism of the space of two-forms:
$W : \Omega^2(M) \rightarrow \Omega^2(M)$.  It will therefore be natural to use the norm associated
to this interpretation, which we denote by using $\| \cdot \|$.  In particular,
$$ |W|^2 = 4 \| W \|^2. $$

\begin{Thme} (Chang-Gursky-Yang \cite{CGY03})
\label{sphere} \newline
Suppose $(M,g)$ is a closed 4-manifold with $g \in \, {\cal A}$.
\noindent \newline
(a). If $\int_M ||W||_g^2 dv_g <  \int _M \sigma _2(A_g)dv_g $
then $M$ is diffeomorphic to either $S^4$ or ${\mathbb {RP}}^4$.\\
%\noindent
(b). If $M$ is not diffeomorphic to $S^4$ or ${\mathbb {RP}}^4$ and
$    \int_M ||W||_g^2 dv_g  =   \int _M \sigma _2(A_g)dv_g $, then 
 $(M, g)$ is conformally equivalent to $({\mathbb {CP}}^2, g_{FS})$.
\end{Thme}

Remark 1: The theorem above is an $L^2$ version of an earlier result of
Margerin \cite{Mar98}. The first part of the theorem should also be compared to
an result of Hamilton \cite{Ha86}; where he pioneered the method of Ricci flow and
established the diffeomorphism of $M^4$ to the 4-sphere under the assumption when the curvature operator is positive.

Remark 2: The assumption $g \in {\cal A}$ excludes out the case when $(M, g) = (S^3 \times S^1, g_{prod})$, where
 $||W||_g = \sigma _2(A_g) \equiv 0$.  
\vskip .1in

\noindent{\underline{ Sketch proof of Theorem \ref{sphere} }}
\vskip .1in 
\begin{proof}
For part (a) of the theorem, we apply the existence argument to find
a conformal metric $g_w$ which satisfies the pointwise inequality
\begin{equation}
 ||W_{g_w} ||^2  < \sigma_2 (A_{g_w})   \,\,\,\,\,  \,\,\,\,  or    \,\,\,\, \, 
\sigma _2(A_{g_w} ) =  ||W_{g_w} ||^2 + c  \,\,\,\text { for  some constant}  \,\,\, c>0.
\end{equation}
The diffeomorphism assertion follows from Margerin's \cite{Mar98} precise
convergence result for the Ricci flow: such a metric will evolve under
the Ricci flow to one with constant curvature. Therefore such a
manifold is diffeomorphic to a quotient of the standard $4$-sphere.\\

For part (b) of the theorem, we argue that if such a manifold is not
diffeomorphic to the 4-sphere, then the conformal structure realizes
the minimum of the quantity $\int |W_g|^2 dv_g$, and hence its Bach tensor
vanishes; i.e. 
$$ B_g = \nabla^k \nabla^l W_{ikjl} + \frac {1}{2} R^{kl} W_{ikjl} =0. $$

As we assume $g \in {\cal  A} $, we can solve the equation
\begin{equation}
\sigma _2(A_{g_w}) =  {(1 - \epsilon) }  ||W_{g_w}||^2 + \, c_{\epsilon},
\end{equation}
where $c_{\epsilon}$ is a constant which tends to zero as $\epsilon$ tends
to zero. We then let $\epsilon$ tends to zero. We obtain in the limit a $C^{1,1}$ metric
which satisfies the equation on the open set $\Omega = \{x| W(x) \neq 0\}$:
\begin{equation}
\label{sharpcp2}
\sigma _2(A_{g_w}) =   ||W_{g_w}||^2.
\end{equation}
We then decompose the $Weyl$ curvature of $g_w$ into its self dual and anti-self dual part as in the Singer-Thorpe decomposition
of the full curvature tensor, apply the
Bach equation to estimate the operator norm of each of these parts as endormorphism on 
curvature tensors, and reduce the 
problem to some rather sophisticated Lagrange multiplier problem. We draw the conclusion that
the curvature
tensor of $g_w$ agrees with that of the Fubini-Study metric
on the open set $\Omega$. Therefore $|W_{g_w}|$ is a constant on $\Omega$,
thus $W$ cannot vanish at all. From this, we conclude that $g_w$ is  Einstein (and under some positive orientation assumption) with $W_{g_w}^-=0$. It follows from a result of Hitchin (see \cite{B87}, chapter 13)
that the limit metric $g_w$ agrees with the Fubini-Study metric of ${\mathbb {CP}}^2$
\end{proof} 
\vskip .1in

We now discuss some of the recent joint work of M. Gursky, Siyi Zhang and myself \cite{CGZ18} extending the theorem \ref{sphere} above
to a perturbation theorem on $\mathbb{CP}^2$.\\

For this purpose, for a metric $g \in {\cal A}$, we define the conformal invariant constant $\beta = \beta ([g]) $ defined as
$$ \int ||W||^2_g dv_g = \beta \int_M \sigma_2 (A_g) dv_g. $$

\begin{lem}
Given $g \in \,{\cal  A} $, if $ 1 < \beta < 2$, then $M^4$ is either homemorphic to either $S^4$ or $\mathbb{RP}^4$ \emph{(hence $b_2^{+} = b_2^{-} =0 $)} or 
$M^4$ is homeomorphic to $\mathbb{CP}^2$ \emph{(hence $b_2^{+} =1 $)}.   
\end{lem}
We remark that $ \beta = 2$ for the product metric on $S^2 \times S^2$.
\\

 %% define b_2+ etc$.
 An additional ingredient to establish the lemma above is the Signature formula :
 $$ 12 \pi^2 \, \tau \,  = \int_M  (||W^{+}||^2 - ||W^{-}||^2) dv ,$$
 where $\tau = b_2^{+} - b_2^{-} $, $||W^{+}||$ is the self dual part of the Weyl curvature and $||W^{-}||$ the anti-self-dual part, $b_2^{+}$, $ b_2^{-}$ the positive and negative part of the intersection form; together
 with an earlier result of M. Gursky 
\cite{Gu98}.\\

 In view of the statement of Theorem \ref{sphere} above, it is tempting to ask if one can change the ``homeomorphism type" 
 to ``diffeomorphism type" in the statement of the Lemma. So far we have not been able to do so, but we have a 
 perturbation result. 

\begin{Thme} (Chang-Gursky-Zhang \cite{CGZ18})
\label{cp2}
There exists some $\epsilon>0$ such that if $(M,g)$ is a four manifold with  $b_2^{+} > 0$ and with a metric of positive Yamabe type satisfying
with $1 < \beta ([g]) <  1 \,+\, \epsilon$, then $(M,g)$ is diffeomorphic to standard $\mathbb{CP}^2$.
 \end{Thme}

\vskip .1in
%\begin{proof}
\noindent {\underline {Sketch proof of Theorem \ref{cp2} }}
\vskip .1in
\begin{proof}
A key ingredient is to apply the condition $b_2^{+}> 0$ to choose a good representative metric $g_{GL} \in [g] $, which is 
constructed in the earlier work of 
Gursky \cite {Gu00} 
and used in Gursky-LeBrun (\cite{GuLeB98}, \cite {GuLeB99}).  To do so, they considered a generalized Yamabe curvature 
%% how to type tilda R $$
$$ {\tilde R }_g = R_g - 2 \sqrt 6 ||W^{+}||_g,$$
and noticed that on manifold of dimension 4, due to the conformal invariance of $||W^{+}||_g$, the corresponding 
Yamabe type functional 
$$
 g \rightarrow   \mu_g :=  \inf_{ g_w\in [g] } \frac {\int_M \tilde R_{g_w} dv_{g_w} }{(vol {g_w})^{\frac {1}{2}} } 
 $$
still attains its infimum; which we denote by $g_{GL}$. The key observation in \cite {Gu00} is that $b_2^{+} >0 $ 
implies $\mu_g <0$ (thus ${\tilde R}_{GL} <0)$.  To see this,  we recall the 
Bochner formula satisfied by the non-trivial self dual harmonic 2-form $\phi$ at the extreme metric:      
$$
\frac{1}{2} \Delta (|\phi|^2) = |\nabla \phi|^2 - 2 W^{+} <\phi, \phi> + \frac{1}{3} R |\phi|^2,
$$  
which together with the algebraic inequality that
$$
- 2 W^{+} <\phi, \phi> + \frac{1}{3} R |\phi|^2 \geq \frac{1}{3} {\tilde R} |\phi|^2,
$$
forces the sign of ${\tilde R}_{GL}$ when $\phi$ is non-trivial. 
\vskip .1in

To continue the proof of the theorem, we notice that for a given metric $g$ satisfying the conformal pinching condition
$1 < \beta ([g])  <  1+ \epsilon$ on its curvature, the corresponding $g_{GL} $ would satisfy $G_2 (g_{GL}) \leq C(\epsilon)$,
where for $k = 2,3,4$, 
$$ 
G_{k} (g) : =  \int_M  \left (  (R- \bar R)^k \, +\,  |Ric^{0} |^k \, + \, ||W^{-}||^k \, + \, |{\tilde R}_{-}|^{k}   \right) dv_g,  
$$ 
where $\bar R$  denotes the average of the scalar curvature  R over the manifold, $Ric^{0} $ the denote the traceless part of the Ricci curvature
and  ${\tilde R}_{-}$ the negative part of $\tilde R$, and where $C(\epsilon)$ is a constant which tends to zero as $\epsilon$ tends to zero.
 
We now finish the proof of the theorem by a contradiction argument and by applying the Ricci flow method of Hamilton
$$ \frac{\partial}{\partial t} g(t) = - 2 Ric_{g(t)} $$
to regularize the metric $g_{GL}$.

Suppose the statement of the theorem is not true, let \{$g_{i} $\} be a sequence of metrics satisfying 
$ 1< \beta ([g_i]) < 1 + {\epsilon_i} $
with $\epsilon_i $ tends to zero as $i$ tends to infinity.  Choose $ g_i(0) = {(g_i)}_{GL} $ to start the Ricci flow, we can derive the inequality
$$
\frac {\partial}{\partial t} G_2 (g(t)) \leq  a G_2( g_t) - b G_4(g(t) ) 
$$ 
for some positive constants $a$ and $b$, also at some fixed time $t_0$ independent of $\epsilon$ and $i$ for each $g_i$.
We then apply the regularity theory of parabolic PDE to derive that some sub-sequences of \{$g_i (t_0)$\} converges to the $g_{FS} $ metric
of $\mathbb{CP}^2$, which in turn implies the original subsequence $(M, \{g_i (0)\})$ hence a subsequence of $(M, \{g_i\})$ is diffeomorphic to $(\mathbb{CP}^2, g_{FS})$.
We thus reach a contradication to our assumption. The reader is referred to the preprint \cite{CGZ18} for details of the proof. 
\end{proof}
\vskip .1in

%%need to check reference
We end this section by pointing out there is a large class of manifolds with metrics in the class ${\cal A}$. 
By the work of Donaldson-Freedman (see \cite{Don83}, \cite{Fr82}) and Lichnerowicz vanishing theorem, the homeomorphism type of the class of simply-connected 4-manifolds which allow a metric with positive scalar curvature 
consists of $S^4$ together with $k \mathbb{CP}^2 \# l {\overline {\mathbb {CP}^2}} $ and $k (S^2 \times S^2)$. Apply some basic algebraic manipulations with the Gauss-Bonnet-Chern formula and the Signature formula, 
we can show that a manifold which admits $g\in\cal{A}$ satisfies $4+5l>k$.
The round metric on $S^4$, the Fubini-Study metric on $\mathbb{CP}^2$, and the product metric on $S^2\times{S^2}$ are clearly in the class $\cal{A}$.
When $l=0$, which implies $k<4$, the class $\cal{A}$ also includes the metrics constructed by Lebrun-Nayatani-Nitta  \cite{LNN97} on $k\mathbb{CP}^2$ for $k\leq2$.
When $l=1$, which implies $k<9$, the class $\cal{A}$ also includes the (positive) Einstein metric constructed by D. Page \cite{P78} on $\mathbb{CP}^2 \# {\overline {\mathbb {CP}^2}} $, the (positive) Einstein metric by Chen-Lebrun-Weber \cite{ChLW08} on $\mathbb{CP}^2 \# 2{\overline {\mathbb {CP}^2}} $,
and the K\"ahler Einstein metrics on $\mathbb{CP}^2 \# l {\overline {\mathbb{CP}^2}}$ for $ 3 \leq l \leq  8 $ as in the work of Tian \cite {T00}. 
It would be an ambitious program to locate the entire class of 4-manifolds with metric in $\cal A $, and to classify their diffeomorphism types by the (relative) size of 
the integral conformal invariants discussed in this lecture.

%We end this section by pointing out there is a large class of manifolds with metrics in the class ${\cal A}$. 
%By the work of Donaldson-Friedman (\cite{Don83}, \cite{Fri82}), the homemorphism type of the class of simply connected 4- manifolds 
%which allows a metric with positive scalar
%curvature consists of $S^4$ together with $k \mathbb{CP}^2 \# l {\overline {\mathbb CP}}^2 $ and $k (S^2 \times S^2)$. 
%Among these manifolds, we know the standard $S^4$, $\mathbb {CP}^2$ and the product metrics on $S^2 \times S^2$ are metrics in class $\cal A$. 
%The list of manifolds which allow metrics in $ {\cal A}$  also includes $ k\mathbb {CP}^2$ for $k \leq 3$ (\cite{ }), $ \mathbb{CP}^2 \# l %
%{\overline{{\mathbb CP}^2 } }$
%for $ l = 1, \,or\, 2$ (\cite{page  }, \cite{  })  as well as  $(CP^2 \#
%l {\overline {\mathbb{CP}^2}}$ for $ 3 \leq l \leq  8 $ which allows positive Kahler Einstein metrics as in the work of Tian \cite {T00}. It would be an ambitious program 
%to find out the entire class of 4-manifolds with allows metric in $\cal A $, and to classify their diffeomorphism types by the (relative) size of its integral conformal invariants.

%%  need to add references for one copy of CP^2 and one copy of its conjugate).

\vskip .2in

\noindent {\bf \S 3. Compact 4-manifold with boundary, $ (Q, T) $ curvatures}
\vskip .1in
\setcounter{section}{3} 
%\begin{document}
%\noindent
%\begin{Thm}
%Suppose $(M,g)$ is a compact oriented manifold of dimension four with positive Yamabe invariant.
%\item{(i)}  If $\int Q_g dv_g = 0$, and if $M$ admits a non-zero harmonic
%\text{\rm 1-}form, then $(M, g)$ is conformal equivalent to a quotient of the
%product space $S^3 \times \mathbb R$. In particular $(M, g)$ is locally
%conformally flat.
%\item{(ii)} If $b_2^+>0$ (i.e. the intersection form has a positive element),
%then with respect to the decompostion of the Weyl tensor into
%the self dual and anti-self dual components $W= W^+ \oplus W^-$,
%\begin{equation}
%\int_M |W_g^+|^2 dv_g \geq \frac{4 \pi ^2}{3}(2 \chi + 3\tau),
%\end{equation}
%where $\tau$ is the signature of $M$. Moreover the equality holds if and only if $g$ is conformal to a
%(positive) Kahler-Einstein
%metric.
%\end{Thm}
\noindent
\setcounter{section}{3}
\setcounter{equation}{00}
\setcounter{Thm}{00}
%\vskip .1in

To further develop the analysis of the $Q$-curvature equation, it is
helpful to consider the associated boundary value problems.
In the case of compact surface with boundary $(X^2, M^1, g)$, where
the metric $g$ is defined on $X^2 \cup M^1$; the Gauss-Bonnet
formula becomes
\begin{equation}
2 \pi \chi(X) = \int_X K \,\,  dv + \oint_M k \,\,  d\sigma,
\end{equation}
where $k$ is the geodesic curvature on $M$. Under conformal
change of metric $g_w$ on $X$, the geodesic curvature changes according to the equation
\begin{equation}
\frac{\partial}{\partial n} w + k \, = \, k_{g_w} \,  e^{w} \,\,\, \text {on M}.
\end{equation}

One can generalize above results to
compact four manifold with boundary $(X^4, M^3 ,g)$; with the role
played by $ (- \Delta,  \frac{\partial}{\partial n})$ replaced by
$(P_4, P_3)$ and with  $(K, k)$ curvature replaced by
$(Q, T)$ curvatures; where $P_4$ is the Paneitz operator and $Q$ the
curvature discussed in section 2; and where $P_3$ is some 3rd order boundary operator 
constructed in Chang-Qing (\cite{CQ97}, \cite{CQ97b}).  The key property of $P_3$ is
that
it is conformally covariant
 of bidegree $(0, 3)$, i.e. 
 $$(P_3)_{g_w} = e^{-3w} (P_3)_{g} $$
when operating on functions defined on the boundary of
 compact $4$-manifolds; and under conformal change of metric
 $ g_w = e^{2w}g$ on $X^4$
we have at the boundary $M^3$
\begin{equation}
P_3 w + T \, = \, T_{g_w} \,  e^{3w}.
\end{equation}
The precise formula of $P_3$ is rather complicated  (see \cite{CQ97}).
Here we will only mention that on
$(B^4, S^3, |dx|^2)$, where $B^4$ is the unit ball in $\mathbb R^4$, we have
\begin{equation}
P_4 = (- \Delta )^2, \,\, P_3 = -\left( \frac{1}{2}\,\,  \frac{\partial}{\partial n} \,\, \Delta + \tilde \Delta \frac{\partial}{\partial n} + \tilde \Delta \right) \,\,\, \,\,\,
\text {and} \,\,\,
T = 2,
\end{equation}
where $\tilde \Delta$ the intrinsic boundary
Laplacian on $S^3$.
In general the formula for $T$ curvature is also lengthy,
$$
T  = \frac 1{12} \frac{\partial}{\partial  n} R + \frac 16 RH - R_{\alpha n \beta
n}L_{\alpha\beta} + \frac 19 H^3 - \frac 13 \text{Tr} L^3 - \frac 13
\tilde \Delta H, 
$$
where $L$ is
the second fundamental form of $M$ in $(X, g)$, and $H$ the mean curvature, and $n $ its the outside normal.  
In terms of these curvatures,  the Gauss-Bonnet-Chern formula can be expressed as:
\begin{equation}
8\pi ^2 \chi (X)= \int _X (||W||^2 \, + \,  2 Q  )\,\, dv \,\,
+ \oint _M ({\cal L} + 2 T ) \,\,  d\sigma.
\end{equation}
where $\cal L$ is a third order boundary curvature invariant that transforms by
scaling under conformal change of metric, i.e. ${\cal L} d \sigma$ is a pointwise conformal invariant.  

The property which is relevant to us is that
$$
\int_X Q dv + \oint_{M} T d\sigma 
$$
is an integral conformal invariant. 

It turns out for the cases which are of interest to us later in this paper, $(X, g)$ is with totally
geodesic boundary, that is, its second fundamental form vanishes. In this special case we have
\begin{equation}
\label{T3}
T = \frac 1{12} \frac{\partial}{\partial  n} R.
\end{equation}
Thus in view of the definitions (2.13)  and (\ref{T3}) of  $Q$ and $T$, in
this case we have
$$
2( \int_X Q dv + \oint_{M} T d\sigma )  = \int_X \sigma_2 \,\, dv, 
$$
which is the key property we will apply later to study the renormalized volume and the compactness problem 
of conformal compact Einstein manifolds in sections 4 and 5 below.
\vskip .2in

{\bf \S 4. {Conformally compact Einstein manifolds}}
\setcounter{section}{4}
\setcounter{equation}{00}
\setcounter{Thm}{00}
\setcounter{lem}{00}
\vskip .1in

{\bf \S 4a. Definition and basics, some short survey}
\vskip .1in

Given a manifold $(M^n, [h]) $, when is it the boundary of a conformally compact Einstein manifold $(X^{n+1}, g^{+})$
with $r^2 g^{+} |_{M} = h$? This problem of finding  ``conformal filling in" is motivated by problems in 
the {AdS/CFT} correspondence in quantum gravity (proposed by Maldacena \cite{Ma98} in 1998)
and from the geometric considerations to study the structure of non-compact asymptotically hyperbolic Einstein manifolds.

Here we will only briefly outline some of the progress made in this problem pertaining to the conformal invariants we are studying.

Suppose that $X^{n+1}$ is a smooth manifold of dimension $n+1$ with
smooth boundary $\partial X = M^n$. A defining function for the
boundary $M^n$ in $X^{n+1}$ is a smooth function $r$ on $\bar X^{n+1}$ such
that
$$
\left\{\aligned r>0 \ \ & \text{in $X$};  \\
r=0 \ \ & \text{on $M$} ;  \\
dr \neq 0 \ \ & \text{on $M$}.
\endaligned\right.
$$
A Riemannian metric $g^{+}$ on $X^{n+1}$ is conformally compact if
$(\bar X^{n+1}, \ r^2 g^{+})$ is a compact Riemannian manifold with
boundary $M^n$ for some defining function $r$. We denote $h := r^2 g^{+} |_{M} $.
\vskip .1in
Conformally compact manifold $(X^{n+1}, g^{+})$ carries a well-defined conformal structure
on the boundary $(M^n, [h])$ by choices of different defining functions $r$. We shall call $(M^n, [h])$ the conformal
infinity of the conformally compact manifold $(X^{n+1}, g^{+})$. 

%A short
%computation yields that, given a defining function $r$,
%$$
%R_{ijkl}[g^+] = |dr|^2_{r^2 g+} ({g+}_{ik}{g+}_{jl}-{g^+}_{il}{g^+}_{jk}) + O(r^3)
%$$
%in a coordinate $(0, \epsilon)\times M^n \subset X^{n+1}$.
%Therefore, if we assume that $g+$ is also asymptotically hyperbolic,
%then
%$$
%|dr|^2_{r^2 g^+}|_M =1
%$$
%for any defining function $r$. 

If $(X^{n+1}, \, g^{+} )$ is a conformally
compact manifold and $\text{Ric}[g^{+}] = -n \, g^{+} $, then we call $(X^{n+1},
\, g^+)$ a conformally compact (Poincare) Einstein (abbreviated as CCE)  manifold. 
We remark that on a CCE manifold $X$, for any given smooth metric $h$ in the conformal infinity $M$, there  
exists a special defining function $r$ (called the geodesic defining function) 
so that 
$
r^2 g^+ |_M = h, \,
$
and $
|dr|^2_{r^2 g^{+}} = 1
$
in a neighborhood of the boundary $[0, \epsilon)\times M$, also the metric $r^2 g^{+}$ has totally geodesic boundary. 
\vskip .1in

\noindent {\underline {Some basic examples}} 
\vskip .1in

{Example 1}:  On $ (B^{n+1}, S^n, g_{\mathbb H}) $
$$
\left(B^{ n+1}, \left( \frac {2}{ 1- |y|^2}\right)^2 |dy|^2\right).
$$
We can then view $(S^n, [g_c]) $ as the compactification of
$B^{n+1}$ using the defining function
$$
 r = 2 \frac {1-|y|}{ 1 + |y|}
$$
$$
g_{\mathbb H} = g^{+} = r^{ -2} \left ( dr^ 2 \, + \, {\left( 1 - \frac
{r^2}{4}\right)}^2 \, g_ c\right).
$$

{Example 2}:  AdS-Schwarzchild space 
$$
On \,\,\, (R^2\times S^2,  g_{m}^+),
$$
where
$$
g_m^{+} = V dt^2 + V^{-1}dr^2 + r^2 g_{c},
$$
$$
V = 1+r^2 - \frac {2m}r,
$$
$m$ is any positive number, $ r\in [r_h, +\infty)$, $t\in
S^1(\lambda)$ and $ g_{c} $ the surface measure on $S^2$ and $r_h$ is the positive
root for $1 + r^2 - \frac {2m}r =0$. 
We remark, it turns out that in this case, there are two different values
of $m$ so that both $g_m^{+} $ are conformal compact Einstein filling for
the same boundary metric $S^1 (\lambda) \times S^2$. This is the famous
non-unique ``filling in" example of Hawking-Page \cite{HaP83}.  
\vskip .1in

\noindent {\underline {Existence and non-existence results}}
\vskip .1in

The most important existence result is the ``Ambient Metric" construction by Fefferman-Graham (\cite{FG85},\cite{FG12}). 
As a consequence of their construction, for any given compact manifold $(M^n, h)$ with an analytic metric $h$, some
CCE metric exists on some tubular neighborhood $M^n \times (0, \epsilon)$ of $M$. 
This later result was recently extended to manifolds $M$ with smooth metrics by Gursky-Sz\'ekelyhidi \cite{GuSz17}.

A perturbation result of Graham-Lee \cite{GL91} asserts that in a smooth neighborhood of the standard 
surface measure $g_c$ on $S^n$, there exist a conformal compact Einstein metric on $B^{n+1}$ with any given conformal infinity $h$.

There is some recent important articles by Gursky-Han and Gursky-Han-Stolz (\cite{GuH17}, \cite{GuHS18}), where they showed that 
when $X$ is spin and of dimension $4k \geq 8$, and the Yamabe invariant $Y(M,[h]) > 0$, then there are topological obstructions 
to the existence of a Poincar\'e-Einstein $g^{+}$ defined in the interior of $X$ with conformal infinity given by $[h]$.  One application of their 
work is that on the round sphere $S^{4k -1}$ with $k \geq 2$, there are infinitely many conformal classes that have no 
Poincar\'e-Einstein filling in in the ball of dimension $4k$.

%The basic idea is to adapt the classical Lichnerowicz result on the vanishing of the $\widehat{A}$-genus for spin manifolds of positive scalar curvature.  
%Indeed, suppose $g_{+}$ is a Poincar\'e-Einstein filling of $[h]$; then one can use the compactification of Lee (see \cite{Lee},\cite{Qing}) to obtain a metric $\overline{g} = \rho^2 g_{+}%$ with positive scalar curvature which is smooth up to the boundary, and such that $M$ is totally geodesic with respect to $\overline{g}$.  It follows that the index of the Dirac operator %(with respect to APS boundary conditions) is zero.  However, using well known properties of the index, it is possible to construct examples of spin manifolds with boundary $M = %
%\partial X$ and conformal classes $[h]$ of positive Yamabe invariant on $M$ such that the index of the Dirac operator (with respect to any extension of any metric in $[h]$) has non-%vanishing index. 

The result of Gursky-Han and Gursky-Han-Stolz was based on a key fact pointed out J. Qing \cite{Q03}, which 
relies on some earlier work of J. Lee \cite{Lee95}.

\begin{lem}
\label{lee}
On a CCE manifold $(X^{n+1}, M^n, g^{+})$, assuming  $Y(M, [h]) >0$, there exists a compactification of $g^{+}$ 
with positive scalar curvature; hence $ Y( X, [r^2 g^{+}]) >0$.  
\end{lem}
\vskip .2in

\noindent{\underline {Uniqueness and non-uniqueness results}}
\vskip .1in

Under the assumption of positive mass theorem, J. Qing \cite{Q03} has established $(B^{n+1}, g_H )$ 
as the unique CCE manifold with $(S^n, [g_c])$ as its conformal infinity. The proof
of this result was later refined and established without using positive mass theorem by Li-Qing-Shi \cite{LiQS17} (see also Dutta and Javaheri \cite{DJ10}).  Later in section 5 of this
lecture notes, we will also prove the uniqueness of the CCE extension of the metrics constructed by Graham-Lee \cite{GL91}
for the special dimension $n=3$.
%% with conformal infinity in a neighborhood of $(S^n, g_c)$ .

As we have mentioned in the example 2 above, when the conformal infinity is $S^1 (\lambda) \times S^2$ 
with product metric, Hawking-Page \cite{HaP83} have constructed non-unique CCE fill-ins. 

\vskip .2in

\noindent {\bf 4b. Renormalized volume}
\vskip .1in

We will now discuss the concept of ``renormalized volume" in the CCE setting, introduced by 
Maldacena \cite{Ma98} (see also the works of Witten \cite{W98}, Henningson-Skenderis \cite{HeK98} and Graham \cite{G00}).
%% add references.
On CCE manifolds $(X^{n+1}, M^n, g^{+})$ with geodesic defining function $r$, \newline 
For $n$ even,
$$ \aligned \text{Vol}_{g^{+}}(\{r
> \epsilon\}) & = c_0 \epsilon^{-n} + c_2\epsilon^{-n+2} + \cdots \\ & +
c_{n-2}\epsilon^{-2} + L \log \frac 1\epsilon + V + o(1).\\
\endaligned $$
For $n$ odd,
$$ \aligned \text{Vol}_{g^{+}}(\{r >
\epsilon\}) & = c_0 \epsilon^{-n} + c_2\epsilon^{-n+2} + \cdots \cdot \cdot  \\
& + c_{n-1} \epsilon^{-1} + V + o(1).\\
\endaligned $$
We call the zero order term $V$ the renormalized volume. It turns out for $n$ even, $L$ is independent of $h \in [h]$ where $h= r^2 g^{+}|_{M}$,  
and for $n$ odd, $V$ is independent of $g \in [g]$, and hence are conformal invariants.

We recall 
\begin{Thm} (Graham-Zworski \cite{GZ03}, Fefferman-Graham \cite{FG02}) \newline 
When n is even, 
$$
L = c_n \oint_M Q_h dv_h,
$$
where $c_n$ is some dimensional constant. 
\end{Thm}

\begin{Thm}(M. Anderson \cite{An01}, Chang-Qing-Yang \cite{CQY06}, \cite{CQY08})
\label{reno} \\
On conformal compact Einstein manifold $(X^4, M^3, g^{+})$, we have
$$ V = \frac {1}{6} \int_{X^4} \sigma_2 (A_g) dv_g $$
for any compactified metric $g$ with totally geodesic boundary. Thus
$$ 
8 \pi^2 \chi (X^4, M^3) = \int ||W||^2_g dv_g + 6 V.
$$
\end{Thm}
\noindent  Remark:  There is a generalization of Theorem \ref{reno} above for any $n$ odd,  with $X^4$ replaced by $X^n$ and with $\int_{X^4} \sigma_2 $ replaced by
some other suitable integral conformal invariants $\int_{X^{n+1}} v^{n+1} $ on any CCE manifold $(X^{n+1}, M^n, g^{+})$; see (\cite {CFGr14}, also \cite{CQY06}).

\vskip .1in

 \noindent{ \underline{Sketch proof of Theorem \ref{reno} for $n=3$}}
 \vskip .1in
\begin{lem} (Fefferman-Graham \cite{FG02})\newline
\label{fgmetric}
Suppose $(X^{4}, M^3, g^{+})$ is conformally compact Einstein with conformal infinity $(M^3, [h])$,  fix $h\in [h]$ and $r$ its
corresponding geodesic defining function.  Consider
\begin{equation}
\label{fg}
 - \Delta_{g^{+}} w \, = \, 3  \,\, \, \, on \,\,\,  X^{4},
\end{equation}
then $ w $ has the asymptotic behavior 
$$ w= log \, r + A + B r^3 $$ near $M$, 
where $A, B$ are functions even in $r$, $A|_M =0$, and
$$ V = \int_M B|_M. $$
\vskip .1in
\end{lem}
 
\begin{lem} 
\label{fgm2} 
With the same notation as in Lemma \ref{fgmetric}, 
consider the metric $g^{*} = g_w = e^{2w} g^+ $, then $g^*$ is totally geodesic on boundary with 
(1)  $Q_{g^{*}} \equiv 0 ,$  
(2)  $ B|_{M} = \frac{1}{36} \frac{\partial}{\partial n} R_{g^*} = \frac {1}{3} T_{g^{*}} .$
\end{lem}
{\underline {Proof of Lemma \ref{fgm2}}:
\begin{proof}
Recall we have $g^+$ is Einstein with $Ric_{g^+}=-3g^+$, thus $$P_{g^+} = (- \Delta_{g^+}) \circ (- \Delta_{g^+} -2) $$
and  $ 2 Q_{g^+} = 6.$
Therefore
$$ P_{g^+} w + 2 Q_{g^+} = 0 = 2 e^{2w} Q_{g^{*}}. $$
Assertion (2) follows from a straight forward computation using the scalar curvature equation and the asymptotic behavior of $w$.
\end{proof} 

Applying Lemmas \ref{fgmetric} and \ref{fgm2}, we get 
$$
\aligned 6 V  & = 6 \oint_{M^3} B|_M d\sigma_{h} 
= \frac {1}{6} \oint_{M^3} \frac{\partial}{\partial n} R_{g^*} d\sigma_h \\
 & = 2 (\int_X Q_{g^{*}} + \oint_M T_{g^{*}} )  = \int_{X^4} \sigma_2 (A_{g^{*}}) dv_{g^{*}}. \\
\endaligned
$$
For any other compactified metric $g$ with totally geodesic boundary, $ \int_{X^4} \sigma_2 (g) dv_g $
is a conformal invariant, and $V$ is a conformal invariant, thus
the result holds once for $g^*$, holds for any such $g$ in the same conformal class, which establishes Theorem \ref{reno}.

\vskip .2in

\noindent{\bf \S 5. {Compactness of conformally compact Einstein manifolds on dimension 4}}
\setcounter{section}{5}
\setcounter{equation}{00}
\setcounter{Thm}{00}
\setcounter{lem}{00}
\setcounter{Prop}{00}
\vskip .1in

In this section, we will report on some joint works of Yuxin Ge and myself {\cite{CGe17}  and also Yuxin Ge, Jie Qing and myself \cite{CGeQ17}}.

The project we work on is to address the problem of given a sequence 
of CCE manifolds $(X^4, M^3, \{(g_i)^{+}\})$ with $M = \partial X$ and $\{ g_i \} = \{ r_i^2 (g_i)^{+} \} $ a sequence of compactified metrics, denote
$h_i = g_i |_M$,  assume $\{h_i\}$ forms a compact family of metrics in $M$, is it true that some representatives ${\bar g_i} \in  [g_i]$
with $\{ {\bar g}_i |_M\} = \{h_i\} $ also forms a compact family of metrics in $X$? 
Let me mention
the eventual goal of the study of the compactness problem is to show existence of conformal filling in for some classes of 
Riemannian manifolds. A plausible 
candidate for the problem to have a positive answer is the class of metrics  $(S^3, h) $ with the scalar curvature of $h$ being positive.
In this case by a result of Marques \cite {Marq12}, the set of such metrics is path-connected, the non-existence 
argument of Gursky-Han, and Gursky-Han-Stolz  (\cite {GuH17}, \cite {GuHS18}) also does not apply. One hopes that our compactness 
argument would lead via either the continuity method or degree theory to the existence of conformal filling in for this class of metrics.
We remark some related program for the problem has been outlined in (\cite{An89}, \cite{An08}).
%We remark that the main goal of the project  is that one hopes the approach
%would lead to existence result of CCE filling in via continuity or degree theory argument. In view of the non-existence result of
%\cite{GuH17} and \cite{GuHS18}, we will target the space $(B^4, S^3)$ with positive scalar curvature on the boundary $S^3$; we
% remark that in this case, 

The first observation is one of the difficulty of the problem is existence of some  ``non-local" term. To see this, we have 
 the asymptotic behavior of the compactified metric $g$ of CCE manifold $(X^{n+1} ,M^n, g^{+})$ with conformal infinity $(M^n, h)$ (\cite{G00}, \cite{FG12}) which in the special case when $n=3$ takes the form
$$  
 g:= r^2 g^{+} = h + g^{(2)} r^2 + g^{(3)} r^3 + g^{(4)} r^4 + \cdot \cdot \cdot \cdot 
$$
on an asymptotic neighborhood of $M \times (0, \epsilon)$, where $r$ denotes the geodesic defining function of $g$. It turns out
$g^{(2)} = - \frac{1}{2}A_{h} $  and is determined by $h$ (we call such terms local terms), $Tr_h g^{(3)} =0 $, while 
$$  
g^{(3)} _{\alpha, \beta} = - \frac {1}{3} {\frac {\partial}{\partial n }{(Ric_{g})}_{\alpha, \beta} } $$
where $\alpha, \beta$ denote the tangential coordinate on $M$, is a non-local term which is not determined by the boundary metric $h$.   
We remark that $h$ together with $g^{(3)}$ determine the asymptotic behavior of $g$ (\cite {FG12}, \cite{Bi08}).

We now observe that different choices of the defining function $r$ give rise to different conformal metric of $\hat h$ in $[h]$ 
on $M$.  For convenience, In the rest of this article, we choose the representative $\hat h = h^{Y}$  be the Yamabe metric with constant scalar curvature in $[h]$ and denote
it by $h$ and its corresponding geodesic defining function by $r$. Similarly one might ask what is a ``good' representative of $\hat g \in [g]$ 
on $X$? Our first attempt is
to choose $\hat g := g^{Y}$, a Yamabe metric in $[g]$. The difficulty of this choice is that it is not clear how to control 
the boundary behavior of $g^{Y}|_M$ in terms of $h^{Y}$. 

We also remark that in seeking the right conditions for the compactness problem, due to the nature of the problem,
 the natural conditions imposed should be conformally invariant.

In the statement of the results below, for a CCE manifold $(X^4, M^3, g^{+})$, and a conformal infinity $(M, [h] )$ with the representative $h = h^{Y} \in [h] $, 
we solve the PDE
\begin{equation}
\label{fg}
- \Delta_{g^{+}} w = 3
\end{equation}
and denote  $g^{*} = e^{2w} g $  be the ``Fefferman-Graham" compactification metric with $g^{*} |_M = h$.

We recall that  $Q_{g^{*}} \equiv 0$, hence the renormalized volume of $(X, M, g^{+}) $ is a multiple of 
$$ 
\int_X \sigma_2 (A_{g*}) dv_{g^{*}} = 2 \oint_M T_{g*} d\sigma_{h} = \frac {1}{6} \oint_M \frac {\partial}{\partial n} R_{g^{*}} d\sigma_{h} .
$$

Before we state our results, we recall formulas for the specific $g^{*}$ metric in a model case.

\begin{lem}
On $(B^4, S^3, g_{\mathbb H})$, 
$$ g^{*} = e^{ (1- |x|^2)}  |dx|^2 \,\,\, on \,\,\, B^4$$  
$$ Q_{g^{*}} \equiv 0, \,\,\,\,\,  T_{g^{*}}  \equiv 2    \,\,\, on \,\,\, S^3 $$
$$ (g^{*})^{(3)} \equiv 0 $$
and
$$ \int_{B^4} \sigma_2 (A_{g*}) dv_{g^{*}}  =  8\, \pi^2. $$
\end{lem} 
 
We will first state a perturbation result for the compactness problem.
 
\begin{Thm} 
\label{maintheorem1}
Let  $\{(B^4, S^3, \{{g_i}^{+}\} )\}$ be a family of oriented CCE on $B^4$ with boundary $S^3$.  
We assume the boundary Yamabe metric $h_i$ in conformal infinity $M$ is of non-negative type. Let $\{g_i^{*}\} $ be the corresponding FG compactification.
Assume  
\begin{enumerate}
\item The boundary Yamabe metrics  $\{h_i\} $ form a compact family in $C^{k+3}$ norm with $k\geq 2$; and  there exists some positive constant $c_1>0$ such that the  Yamabe constant for the conformal infinity $[h_i]$ is bounded uniformly from below by $c_1$, that is, 
$$
Y(M,[h_i])   \ge c_1;
$$
\item There exists some small positive constant $\varepsilon>0$ such that for all $i$
\begin{equation}
\label{s2}
\int_{B^4} \sigma_2 (A_{g_i^*}) dv_{g_i^{*}} \, \geq  \,  8 \pi^2 - \varepsilon. 
\end{equation} 
\end{enumerate}
Then the family of the $g_i^{*} $  is compact in $C^{k+2,\alpha}$ norm for any $\alpha\in (0,1)$ up to a diffeomorphism fixing the boundary.
\end{Thm}

Before we sketch the proof of the theorem, we will first mention that on $(B^4, S^3, g )$, 
for a compact metric $g$  with totally geodesic boundary, Gauss-Bonnet-Chern
formula takes the form:
$$
8 \pi^2 \chi (B^4, S^3) = 8 \pi^2 = \int_{B^4} ( ||W||_g^2 + \sigma_2 (A_g)) dv_g,  
$$
which
 together with the conformal invariance of the $L^2$ norm of the Weyl tensor,  
imply that the condition (\ref{s2}) in the statement of the Theorem \ref{maintheorem1} is equivalent to
\begin{equation} 
\label{w2}
\int_{B^4} ||W||_{g_i^{+}} ^2 dv_{g_i^{+}} \leq \varepsilon.
\end{equation}

What is less obvious is that in this setting, we also have other equivalence conditions as stated in Corollary (\ref{application1}) below. 
This is mainly due to following result by Li-Qing-Shi (\cite{LiQS17}).
\begin{Prop}
\label{collapse} 
Assume that $(X^{n+1}, g^+)$ is a CCE manifold with $C^3$ regularity whose conformal infinity is of positive Yamabe type. Let $p \in X$ be a fixed point and $t>0$. Then
$$
\left(\frac{  Y(\partial  X,[h])}{Y(\mathbb{S}^{n}, [g_c])}\right)^{\frac{n}{2}}\le  \frac{Vol(\partial B_{g^{+}} (p,t))}{ Vol(\partial B_{g_\mathbb{H}} (p,t))}\le  \frac{Vol( B_{g^{+}} (p,t))}{ Vol( B_{g_\mathbb{H}} (p,t))}\le 1
$$
where $B_{g^{+}} (p,t)$ and $B_{g_\mathbb{H}} (p,t)$ are geodesic balls.
\end{Prop}

\begin{cor} 
\label{application1}
Let  $\{X=B^4,M=\partial X=S^3, g^+\}$ be a  $4$-dimensional oriented CCE on $X$ with boundary $\partial X$.  
Assume the boundary Yamabe metric $h= h^Y$ in the conformal infinity of positive type and $Y(S^3, [h])> c_1$ for some fixed $c_1>0$ and  $h$ is bounded 
in $C^{k+3}$ norm with $k\geq 5$. Let $g^*$ be the corresponding FG compactification.
Then the following properties are equivalent: 
\begin{enumerate}
\item There exists  some small positive number $\varepsilon>0$ such that
\begin{equation}
\int_X \sigma(A_{g*}) dv_{g^{*}}  \geq 8\pi^2- \varepsilon.
\end{equation} 
\item There exists some small positive  number $\varepsilon>0$ such that
$$
\int_X ||W||_{g^+} ^2 dv_{g^+}  \le \varepsilon .
$$
\item There exists  some small positive  number $\varepsilon_1>0$ such that
$$
Y(S^3,[g_{c}])\, \geq \,  Y(S^3,[h])> Y(S^3,[g_{c}])- \varepsilon_1
$$
where $g_{c}$ is the standard metric on $S^3$.
\item There exists some small positive number $\varepsilon_2>0$ such that for all metrics $g^{*}$ with boundary metric $h$ same
volume as the standard metric $g_c$ on $S^3$, we have
$$ 
T( g^{*}) \geq 2 - \varepsilon_2.
$$
\item There exists  some small positive  number $\varepsilon_3>0$ such that
$$
|(g^{*})^{(3)}| \leq \varepsilon_3.
$$
\end{enumerate}
Where all the $\varepsilon_{i} $ (i = 1,2,3) tends to zero when $\varepsilon$ tends to zero and vice versa for each $i$. 
%\end{enumerate}
\end{cor}

Another consequence of Theorem \ref{maintheorem1} is  the ``uniqueness" of the Graham-Lee metrics mentioned in section 4a.

\begin{cor}
\label{maincorollary0}
There exists some $\varepsilon>0$, such that for all metrics $h$ on $S^3$ with  $ || h - g_c||_{C^{\infty}}  < \varepsilon $, 
there exists a {\it unique} CCE filling in $(B^4, S^3, g^{+} )$ of $h$.   
\end{cor}

\noindent {\underline {Sketch proof of Theorem \ref{maintheorem1}}}
\vskip .1in
We refer the readers to the articles Chang-Ge \cite{CGe17} and Chang-Ge-Qing \cite{CGeQ17}, both will soon be posted on arXiv
for details of the arguments, here we will present a brief outline.  
\vskip .1in

We first state a lemma summarizing some analytic properties of the metrics $g^{*}.$
\begin{lem}
\label{fgcom}
On a CCE manifold $(X^4, M^3, g^+)$, where the scalar curvature of the
conformal infinity $(M, h)$ is positive. Assume $h$ is at least $C^l$
smooth for $l \geq 3$.  Denote $g^{*} = e^{2w} g^+ $ the FG compactification. Then \\
(1) $Q (g^{*}) \equiv 0$, \newline
(2) $ R_{g^{*}} >0$, which implies in particular $ |\nabla_{g^{*}} w | e^{w}
\leq 1$. \\
(3) $g^{*}$ is Bach flat and satisfies an $\epsilon$-regularity property,
which implies in particular, once it is $C^3$ smooth, it is
$C^{l}$ smooth for $ l \geq 3$.
\end{lem}

We remark that statement (2) in the Lemma above follows from 
a continuity argument via some theory of scattering matrix (see Case-Chang \cite {CaC16}, \cite{CaC16b}), with the starting point of the 
argument the positive scalar curvature metric constructed by J. Lee which we have mentioned earlier in Lemma \ref{lee}.  \\

{\underline {Sketch proof of Theorem \ref{maintheorem1}} \\

Proof of the theorem is built on contradiction arguments. We first note, assuming the conclusion of the theorem does not hold, then there is a sequence of \{$g_i^*$\}
which is not compact so that the $L^2$ norm of its Weyl tensor of the sequence tends to zero.\\ 

Our main assertion of the proof is that the $C^1$ norm of the curvature of the family \{${g_i }^{*} $\}  remains uniformly bounded. \\

Assume the assertion is not true, we rescale the metric
$
\bar g_i= K_i^2 {g_i}^{*}
$
where there exists some point $p_i\in X$ such that
$$
K_i^2=\max\{\sup_{B^4}  |Rm_{{g_i}^{*}} |, \sqrt{\sup_{B^4} |\nabla Rm_{{g_i}^{*} } |}\}=
|Rm_{{g_i}^{*}} |(p_i) (\mbox{ or } \sqrt{ |\nabla Rm_{{g_i}^{*}} |(p_i)})
$$
We mark the accumulation point $p_i$ as $0\in B^4$. Thus, we have
\begin{equation}
\label{marking}
|Rm_{\bar g_i}|(0)=1 \mbox{ or } |\nabla Rm_{\bar g_i}|(0)=1
\end{equation}
We denote the corresponding defining function $\bar w_i$ so that 
$\bar g_i = e^{ 2 \bar w_i} g^+_i$; that is $e^{ 2 \bar w_i} = K_i^2 e^{2 w_i}$  and denote
$\bar h_i := \bar g_i|_{S^3}$. We remark that the metrics $\bar g_i$ also satisfy the conditions in Lemma \ref{fgcom}.  

As $0\in \,  {\overline {B^4}} $  is an accumulation point, depending on the location of $0$, we call it either an interior or a boundary blow up
point; in each case, we need a separate argument but for simplicity here we will assume we have $0\in \,  S^3$ is a boundary blow up point. In this case, we  
denote the $ (X_{\infty} , g_{\infty}) $ the Gromov-Hausdorff limit of the sequence $(B^4, \bar g_i) $, and $h_{\infty} := g_{\infty}|_{S^3}$. \\

Our first observation is that it follows from the assumption (1) in the statement of the theorem, we have $( \partial X_{\infty}, h_{\infty}) = (\mathbb {R}^3, |dx|^2)$. \\

Our second assertion is that by the estimates in Lemma \ref{fgcom},  one can show $\bar w_i$ converges uniformly on compacta
on $X_{\infty}$, we call the limiting function $ \bar w_{\infty}$. Hence, the corresponding metric 
$ g_{\infty}^{+} $ with $g_{\infty} = e^{2 {\bar w}_{\infty}} g_{\infty}^+$ exists, satisfying $Ric_{ g_{\infty}^{+}} = - 3 g_{\infty}^{+}$,  i.e., 
the resulting $g_{\infty}$ is again conformal to a Poincare Einstein metric 
$g_{\infty}^{+}$ with $ ||W_{g_{\infty}^+} || \equiv  0$. \\

Our third assertion is that  $ (X_{\infty},  g_{\infty}^{+}) $ is (up to an isometry) the model space $( \mathbb {R}^4_{+}, \, g_{\mathbb H} := \frac {|dx|^2 + 
|dy| ^2}{ y^2}) $, where $ \mathbb {R}^4_{+} = \{ (x, y) \in \mathbb {R}^4 | y > 0 \}$. We can then 
apply a Liouville type PDE argument to conclude $\bar  w_{\infty} = \log \,  y $. \\

Thus $g_{\infty}$ is in fact the flat metric $ |dx|^2 + |dy|^2$, which contradicts the marking property (\ref{marking}). 
This contradiction establishes the main assertion.\\

Once the $C^1$ norm of the curvature of the metric \{$g_i^{*}$\}  is bounded, we can apply some further blow-up argument to show the 
diameter of the sequence of metrics is
uniformly bounded, and apply a version of the Gromov-Hausdorff compactness result (see \cite{ChGrT82})
 for compact manifolds with totally geodesic boundary to prove that \{$g_i^{*} $\}  forms a compact family in a suitable $C^l$ norm for some $ l\geq 3$. This finishes the proof of the theorem.

%\end{proof}

\vskip .1in

We end this discussion by mentioning that, from Theorem \ref{maintheorem1}, an obvious question to ask is if we can extend the
perturbation result by improving condition (\ref{s2}) to the condition $ \int_{B^4}  \sigma_2 >0$. This is a direction we are working on but 
have not yet be able to accomplish. Below are statements of two more general theorems that we have obtained in (\cite{CGe17}, \cite{CGeQ17}). 

 \begin{Thm}
\label{maintheorem2}
Under the assumption (1) as in Theorem \ref{maintheorem1}, assume further 
the $T$ curvature on the boundary $T_i=  \frac{1}{12} \frac{\partial} {\partial n}  R_{{g_i}^{*}}$  satisfies the following condition
\begin{equation}
\label{TT3} 
\liminf_{r\to 0}\inf_i\inf_{S^3}\oint_{\partial B(x,r)}T_i\ge 0.
\end{equation}
Then the family of metrics \{$g_i^{*} $\}  is compact in $C^{k+2,\alpha}$ norm for any $\alpha\in (0,1)$ up to a diffeomorphism fixing the boundary, provided $k\ge 5$.
\end{Thm}

\begin{Thm}
\label{maintheorem3}
Under the assumption (1)  as in Theorem \ref{maintheorem1}, assume further that  
there is no concentration of $ S_i  := \, (g_i^{*})^{(3)}  $-tensor defined on $S^3$  in $L^1$ norm for the $g_i^{*} $ metric
in the following sense,
\begin{equation} 
\label{Stensor}
\lim_{r\to 0}\sup_i\sup_x\oint_{\partial B(x,r)} |S_i|=0.
\end{equation}
Then, the family of the metrics \{$g_i^{*} $\}   is compact in $C^{k+2,\alpha}$ norm for any $\alpha\in (0,1)$ up to a diffeomorphism fixing the boundary,
provided $k \geq 2$.
\end{Thm}

The reason we can pass the information from the $T$ curvature in Theorem \ref{maintheorem2} to the $S$ tensor in Theorem \ref{maintheorem3} is due to the fact that 
for 
the blow-up limiting metric $g_{\infty}$, $T_{g_{\infty}} \equiv 0$ if and only if $ S_{g_{\infty}} \equiv 0$. \\
%The later fact was established 
% earlier in Case-Chang (\cite{CaC16}, \cite{CaC16b}) based on some argument  on CCE manifolds.

It remains to show if there is a connection between condition (\ref{TT3}) in Theorem \ref{maintheorem2} to the positivity of the renormalized volume, i.e. when
$\int_X \sigma_2(A_g) dv_g > 0$.

%\end{document}

\newcommand{\namelistlabel}[1] {\mbox{#1}\hfil}
\newenvironment{namelist}[1]{%
\begin{list}{}
{\let\makelabel\namelistlabel
\settowidth{\labelwidth}{#1}
\setlength{\leftmargin}{1.1\labelwidth}}
}{%
\end{list}}

%\begin{document}

%\newpage

%\end{thebibliography}
\vskip .3in
\noindent
{\bf Address:}
\noindent
\begin{namelist}{xxxxxxxxxxxxxxxxxxxxxxxxx}
\item[{\sc Sun-Yung Alice Chang,}]
Department of Mathematics, \\
Princeton University, \\
Princeton, NJ 08544  \\
{\it email: chang@math.princeton.edu}

\end{namelist}


\small

\begin{thebibliography}{99}



%\bibitem{AES18} K. Akutagawa,  H.  Endo and  H. Seshadri; {\em  A gap theorem for positive Einstein metrics on the four-sphere}, preprint Arxiv: 1801.10305V2.

\bibitem{An01} M. Anderson; {\em $L^2$ curvature and volume renormalization of
the AHE metrics on 4-manifolds}, Math. Res. Lett., \textbf{8} (2001) 171-188.

\bibitem{An89} M. Anderson, {\em Ricci curvature bounds and Einstein metrics  on  compact manifolds}, J. Amer. Math. Soc., \textbf{2} (1989), no. 3, 455-490.

\bibitem{An08} M. Anderson; {\em Einstein metrics with prescribed conformal infinity on 4-manifolds}, Geom. Funct. Anal., \textbf{18} (2008) 305-366.

%\bibitem{An89} M. Anderson, {\em Ricci curvature bounds and Einstein metrics  on  compact manifolds}, JAMS (1989) 455-490.

%\bibitem{A} D. Adams; {\em A sharp inequality of J. Moser for
%higher order derivatives}, Ann. of Math., 128 (1988), 385-398.

%\bibitem{Al} O. Alvarez; {\em Theory of strings with boundary},
%Nucl. Phys.,  B 216 (1983), 125-184.

\bibitem{Au76} T. Aubin; {\em Equations differentielles non lineaires et
probleme de Yamabe concernant la courbure scalaire}, J. Math. Pures
Appl. \textbf{55} (1976), 269-296.

%\bibitem{Au79} T. Aubin; {\em Meilleures constantes dans le theorem
%d'inclusion
%de Sobolev et un theorem de Fredholme non lineaire pour la transformation
%conforme de la courbure scalaire}, J. Funct. Anal., 32 (1979), 148-174.

%\bibitem{BC} A. Bahri and J. M. Coron; {\em The Scalar curvature
% problem on the standard three dimensional sphere}, 
% J. Funct. Anal., 95 (1991), 106-172.

%\bibitem{Be93} W. Beckner; {\em Sharp Sobolev inequalities on the sphere
% and the Moser-Trudinger inequality}, Ann. of Math., 138 (1993), 213-242.
\bibitem{B87} A. Besse; {\em Einstein Manifolds}, Berlin: Springer-Verlag (1987).

%\bibitem{Biquard1} O. Biquard, {\sl Einstein deformations of hyperbolic metrics}

\bibitem{Bi08} O. Biquard; {\em Continuation unique \`a partir de l'infini conforme pour les m\'etriques d'Einstein}, Math. Res. Lett., \textbf{15} (2008), 1091-1099.


%\bibitem{BH14} O. Biquard and M. Herzlich; {\em  Analyse sur un demi-espace hyperbolique et poly-homog\'en\'eit\'e locale}, Calc. Var. P.D.E., 51 % (2014), 813-848.

%O. Biquard and M. Herzlich; {\em  Analyse sur un demi-espace hyperbolique et poly-homog\'en\'eit\'e locale}, \newblock {\em Calc. Var. P.D.E. 51} (2014), 813-848.

\bibitem{Br85} T. Branson; {\em Differential operators canonically
associated to a conformal structure}, Math. Scand., \textbf{57} (1985), 293-345.


\bibitem{Br95} T. Branson; {\em Sharp inequality, the functional
determinant and the complementary series}, Trans. Amer. Math. Soc., \textbf{347} (1995), 3671-3742.

\bibitem{Br96} T. Branson; {\em An anomaly associated with
4-dimensional
quantum gravity}, Comm. Math. Phys., \textbf{178} (1996), 301-309.

%\bibitem{Br-5} T. Branson; {\em The functional determinant},
% Lecture notes series, No. 4, Seoul National University.

\bibitem{BCY92} T. Branson, S.-Y. A. Chang and P. Yang; {\em Estimates
and extremals for the zeta-functional determinants on four-manifolds},
Comm. Math. Phys., \textbf{149} (1992), no. 2, 241-262.

\bibitem{BO91} T. Branson and B. {\O}rsted; {\em Explicit functional
determinants in four dimensions}, Proc. Amer. Math. Soc., \textbf{113}
(1991), 669-682.

%\bibitem{Bre02} S. Brendle; {\em Global existence and convergence for a
%higher order
%flow in conformal geometry}, preprint, 2002.


\bibitem{BV04} S. Brendle and J. Viaclovsky; {\em A variational characterization for $\sigma_{n/2}$}, Calc. Var. Partial Differential
Equations, \textbf{20} (2004), no. 4, 399-402.


%\bibitem{Ca90} L. A. Caffarelli; {\em A localization property of
%viscosity
%solutions to the Monge-Amp\'ere equation and their strict
%convexity},  Ann. of Math. (2) 131 (1990), no. 1, 129-134.


%\bibitem{CNS85} L. Caffarelli, L. Nirenberg and J. Spruck; {\em The
%Dirichlet problem for nonlinear second order elliptic equations, III:
%Functions of the eigenvalues of the Hessian}, Acta Math., 155
%(1985), no. 3-4, 261-301.

%\bibitem{CaYa95} L. Caffarelli and Y. Yang; {\em Vortex
%condensation in the
%Chern-Simons Higgs model: an existence theorem},
%Comm. Math. Phys., 168 (1995), 321-336

%\bibitem{CC86} L. Carleson and S.-Y. A. Chang; {\em On the existence
%of an extremal
%function for an inequality of J. Moser}, Bull. Sci. Math.,
%110 (1986), 113-127.
 
\bibitem{CaC16} J. Case and S.-Y. A. Chang;  {\em On fractional GJMS operators},  Comm. Pure Appl. Math., \textbf{69} (2016), 1017-1061.

\bibitem{CaC16b} J. Case and S.-Y. A. Chang; {\em  Errata for ``On fractional GJMS operators"}, arXiv:1406.1846.

\bibitem{CaY13} J. S. Case and P. C. Yang; {\em A Paneitz-type operator for CR pluriharmonic functions}, Bull. Inst. Math. Acad. Sin. (N.S.) \textbf{8} (2013), no. 3, 285-322.

\bibitem{CGe17} S.-Y. A. Chang and Y. Ge; {\em  Compactness of conformally compact Einstein manifolds in dimension 4}, preprint.

\bibitem{CGeQ17} S.-Y. A. Chang, Y. Ge, and J. Qing; {\em  Compactness of conformally compact Einstein manifolds in dimension 4, part II}, preprint.

\bibitem{CFGr14} S.-Y. A. Chang, H. Fang, and C. R. Graham;  {\em A note on renormalized volume functionals}, Differential Geometry and its Applications, 33, Supplement, March 2014, 246-258.
Special issue in honor of Michael Eastwood's 60th birthday.

%\bibitem{CGY93} S.-Y. A. Chang, M. Gursky and P. Yang; {\em Prescribing
%scalar curvature on $S^2$ and $S^3$}, Calculus of Variation, 1 (1993),
 %205-229.

%\bibitem{CGY99}  S.-Y. A. Chang, M. Gursky and P. Yang; {\em
%Regularity of a fourth order PDE
%with critical exponent}, Amer. J. of Math.,  121 (1999)  215-257.

\bibitem{CGY02} S.-Y. A. Chang, M. Gursky, and P. Yang; {\em An equation
of Monge-Amp\`ere type in conformal geometry, and four-manifolds of positive Ricci
curvature},  Ann. Math., \textbf{155} (2002), no. 3, 711-789.

%\bibitem{CGY02b}  S.-Y. A. Chang, M. Gursky and P. Yang; {\em An a prior
%estimate for a fully nonlinear equation on Four-manifolds},
%J. D'Analyse Math.; Thomas Wolff memorial issue, 87 (2002), .

%\bibitem{CGY}  S.-Y. A. Chang, M. Gursky and P. Yang; {\em Entire solutions of a
%fully nonlinear equation}, 


\bibitem{CGY03} S.-Y. A. Chang, M. J. Gursky, and P. Yang; { \em A conformally invariant sphere theorem in four dimensions}, 
Publ. Math. Inst. Hautes \'Etudes Sci. \textbf{98} (2003), 105-143.

\bibitem{CGZ18} S.-Y. A. Chang, M. J. Gursky, and S. Zhang; {\em Conformal gap theorems on $\mathbb {CP}^2$}, preprint.  

\bibitem{CQ97} S.-Y. A. Chang and J. Qing; {\em The Zeta functional
determinants on manifolds with boundary I--the formula}, J. Funct. Anal.,
\textbf{147} (1997), no.2, 327-362.

\bibitem{CQ97b} S.-Y. A. Chang and J. Qing; {\em The Zeta functional
determinants on manifolds with boundary II--Extremum metrics and
compactness of isospectral set}, J. Funct. Anal., \textbf{147} (1997), no.2, 363-399.

%\bibitem{CQY00}  S.-Y. A. Chang, J. Qing and P. Yang; {\em On the
%Chern-Gauss-Bonnet integral for conformal metrics on $R^4$}, Duke Math J,
%103, No. 3 (2000), 523-544.

\bibitem{CQY06}
S.-Y. A. Chang, J. Qing, and P. Yang; {\em On the renormalized volumes for
conformally compact Einstein manifolds},  (Russian) Sovrem. Mat. Fundam. Napravl. \textbf{17} (2006), 129-142.

\bibitem{CQY08}
S.-Y. A. Chang, J. Qing, and P. Yang;  {\em   On the renormalized volumes for conformally
  compact einstein manifolds}, J. Math. Sci. \textbf{149}  (2008), 1755-1769.


%\bibitem{CQY-2}  S.-Y. A. Chang, J. Qing and P. Yang; {\em Compactification
%for a class of conformally flat 4-manifold},  Inventiones
%Mathemticae, 142 (2000), 65-93.

% \bibitem{CQY04}
%S.Y. A. Chang, J. Qing and P. Yang, {\em On the topology of conformally
%compact Einstein 4-manifolds}, \newblock Noncompact Problems at the intersection
%of Geometry, Analysis and Topology \newblock {\em Contemporary Math.} 350 (2004),  49-61.

%\bibitem{CQY04b}
%S.Y. A. Chang, J. Qing and P. Yang,  {\em   On the topology of conformally compact Einstein 4-manifolds}, Noncompact problems at the intersection of geometry, analysis, and topology, %Contemp. Math., vol. 350, Amer. Math. Soc., Providence, RI, 2004, pp. 49-61. MR 2082390, https://doi.org/10.1090/conm/350/06337



%\bibitem{CLW17}
%X. Chen, M. Lai   and F. Wang; {\em  Escobar-Yamabe compactifications for Poincar\'e-Einstein manifolds and rigidity theorems},
%preprint arxiv 1712.02540 (2017).


\bibitem{CY95} S.-Y. A. Chang and P. Yang; {\em Extremal metrics of
zeta functional determinants on 4-Manifolds}, Ann. Math.,
\textbf{142} (1995), 171-212.

%\bibitem{CY97}  S.-Y. A. Chang and P. Yang; {\em On uniqueness of solution of a n-th order differential equation in conformal geometry}, Math. Res
%Letter, 4 (1997), 91-102.

\bibitem{CY03} S.-Y. A. Chang and P. Yang; {\em The inequality of Moser
and Trudinger and applications to conformal geometry}, CPAM, Vol LVI, 
August (2003), no. 8, 1135-1150. Special issue dedicated to the memory of
J\"urgen K. Moser.

%Comm. Pure Appl. Math., special issue in memory of J\"urgen Moser.

\bibitem{CQY07} S.-Y. A. Chang, J. Qing, and P. Yang; {\em On a conformal gap and finiteness theorem for a class of four-manifolds}, Geom. Func. Anal. \textbf{17} (2007), 404-434.

\bibitem{Ch70} J. Cheeger;  \newblock {\em Finiteness theorem for riemannian manifolds}, Amer. J. Math., \textbf{92} (1970), 61-74.

\bibitem{ChGrT82} J. Cheeger, M. Gromov, and M. Taylor; {\em Finite propagation speed, kernel estimates for functions of the Laplace operator, and the geometry of complete Riemannian manifolds}, J. Diff. Geom., \textbf{17} (1982), 15-53.

\bibitem{ChNa15} J. Cheeger and A. Naber; \newblock {\em Regularity of Einstein manifolds and the codimension 4 conjecture}, Ann. Math., \textbf{182} (2015), 1093-1165.

%\bibitem{Chen05b} S. Chen; \newblock {\em Conformal deformation on manifolds with boundary}, \newblock {Geom. Funct. Anal.} 19 (2009), 1029-1064.


%\bibitem{ChLW17}
%X. Chen, M. Lai   and F. Wang; {\em Escobar-Yamabe compactifications for Poincar\'e-Einstein manifolds and rigidity theorems},
% preprint arxiv 1712.02540 (2017).

\bibitem{ChLW08} X. Chen, C. Lebrun, and B. Weber; {\em On conformally K\"ahler, Einstein manifolds}, J. Amer. Math. Soc., \textbf{21} (2008), 1137-1168.

%\bibitem{DM08}
%D.Djadli and A. Malchiodi;
 %\newblock {\em Existence of conformal metrics with constant $Q$-curvature}, \newblock { Ann. of Math.} 168 (2008), 813-858.

\bibitem{DJ10} S. Dutta and M. Javaheri; \newblock {\em Rigidity of conformally compact manifolds with the round sphere as the conformal infinity}, Adv. Math., \textbf{224} (2010), 525-538

\bibitem{Don83} S. Donaldson; {\em An application of gauge theory to four-dimensional topology}, J. Diff. Geom., \textbf{18} (1983), 279-315.

\bibitem{DK90} S. Donaldson and P. Kronheimer; {\em The Geometry of Four-manifolds}, 
Oxford Univ. Press, New York, 1990.

\bibitem{ES85} M. Eastwood and M. Singer; {\em A conformally invariant
Maxwell gauge}, Phys. Lett. A, \textbf{107} (1985), 73-83.

%\bibitem{Ev82} C. Evans; {\em Classical solutions of fully non-linear,
%convex, second order elliptic equations},  Comm. Pure Appl. Math.,
 %XXV (1982), 333-363.


\bibitem{FG85} C. Fefferman and C. R. Graham; {\em Conformal
invariants},
In: {\it \'{E}lie Cartan et les Math\'{e}matiques d'aujourd'hui.
Asterisque} (1985), 95-116.



\bibitem{FG02} C. Fefferman and C. R. Graham;
\newblock {\em  $Q$-curvature and Poincar\'e metrics}, \newblock {Math. Res. Lett.,} \textbf{9} (2002), 139-151.

\bibitem{FG12} C. Fefferman and C. R. Graham;
\newblock {\em  The ambient metric},  Annals of Mathematics Studies, 178, Princeton University Press, Princeton, (2012).


%\bibitem{GiPR79} G. Gibbons, C. Cope and A. Romer  {\sl  Index theorem boundary terms for gravitational instantons}, \newblock { Nuclear Phys. }B157 (1979), 377-386.

\bibitem{Ga59} L. G{\aa}rding; {\em An inequality for hyperbolic polynomials},
 J. Math. Mech., \textbf{8} (1959), 957-965.

\bibitem{Fr82} M. Freedman; {\em The topology of four-dimensional manifolds}, J. Diff. Geom., \textbf{17} (1982), 357-453.

\bibitem{G00} C. R. Graham;
\newblock {\em  Volume and Area renormalizations for conformally compact Einstein metrics}, \newblock The Proceedings of the 19th Winter School ``Geometry and Physics" (Srn\`{i}, 1999).
\newblock Rend. Circ. Mat. Palermo \textbf{63} (2000), 31-42.

\bibitem{GL91} C. R. Graham and J. Lee;
\newblock {\em  Einstein metrics with prescribed conformal infinity on the ball},  \newblock { Adv. Math.,} \textbf{87} (1991), 186-225.

\bibitem{GJMS92} C. R. Graham, R. Jenne, L. Mason, and G. Sparling;
{\em Conformally invariant powers of the Laplacian, I: existence},
J. London Math. Soc., \textbf{46} (1992), no. 2, 557-565.

\bibitem{GZ03} C. R. Graham and M. Zworski;
\newblock {\em  Scattering matrix in conformal geometry}, \newblock {Invent. Math.,} \textbf{152} (2003), 89-118.

\bibitem{Gu98} M. J. Gursky; {\em The Weyl functional, de Rham cohomology, and K\"ahler-Einstein metrics}, Ann. Math., \textbf{148} (1998), 315-337.

\bibitem{Gu99} M. J. Gursky; {\em The principal eigenvalue of a conformally invariant differential operator, with an Application to semilinear elliptic PDE}, Comm. Math. Phys., \textbf{207} (1999), 131-143.

\bibitem{Gu00} M. J. Gursky;  {\em Four-manifolds with $\delta{W}^+=0$ and Einstein constants of the sphere}, Math. Ann., \textbf{318} (2000), 417-431.

\bibitem{GuH17} M. J. Gursky  and Q. Han; {\em Non-existence of Poincare-Einstein manifolds with prescribed conformal infinity},  Geom. Funct. Anal., \textbf{27} (2017), no. 4, 863-879.

\bibitem{GuHS18} M. J. Gursky, Q. Han, and  S. Stolz; {\em  An invariant related to the existence of conformally compact Einstein fillings}, preprint.

\bibitem{GuHL16} M.J. Gursky, F. Hang, and Y.-J. Lin; {\em Riemannian Manifolds with Positive Yamabe Invariant and Paneitz Operator}, 
Int. Math. Res. Not., (2016), no.5, 1348-1367. 


\bibitem{GuLeB98} M.J. Gursky and C. Lebrun; {\em Yamabe invariants and
Spinc structures},  Geom. Funct. Anal., \textbf{8} (1988), no. 6, 965-977.

\bibitem{GuLeB99} M. J. Gursky and C. Lebrun; {\em On Einstein manifolds of positive sectional curvature}, Ann. Glob. Anal. Geom. \textbf{17} (1999), 315-328.


\bibitem{GuM15}
M. J. Gursky and A. Malchiodi; \newblock {\em   A strong maximum principle for the {P}aneitz
  operator and a non-local flow for the {$Q$}-curvature}, \newblock { J. Eur. Math. Soc.,} \textbf{17} (2015), 2137-2173.

\bibitem{GuSt16} M.J. Gursky and J. Streets; {\em A formal Riemannian structure on conformal classes and uniqueness for the $\sigma_2$ -Yamabe
problem}, arkiv:1603.07005v1.

\bibitem{GuSz17} M. J Gursky and G. Szekelyhidi;  {\em A local existence result for
Poincar\'e-Einstein metrics}, arXiv:1712.04017. 


\bibitem{GV01} M. Gursky and J. Viaclovsky; {\em A new variational
characterization of three-dimensional space forms},
Invent. Math., \textbf{145} (2001), 251-278.

\bibitem{Ha86} R. Hamilton; {\em Four manifolds with positive curvature
operator}, J. Diff. Geom.,  \textbf{24} (1986), 153-179.

\bibitem{Ham82} R. Hamilton; {\em Three-manifolds with positive Ricci curvature}, J. Diff. Geom., \textbf{17} (1982), 255-306.

\bibitem{HY15} F. Hang and P. C. Yang; {\em Sign of Green's function of Paneitz operators and the Q-curvature}, Int. Math. Res. Not., \textbf{19} (2015), 9775-9791.

\bibitem{HY16} F. Hang and P. C. Yang; {\em Lectures on the fourth-order Q-curvature
equation} in {\em Geometric Analysis around Scalar Curvatures}, Lect. Notes
Ser. Inst. Math. Sci. Natl. Univ. Singap., World Sci. Publ., Hackensack, NJ, \textbf{31} (2016), 1-33.

%\bibitem{HY16b} F.  Hang and P. C.�Yang; {\em�Q�curvature on a class of 3-manifolds}, �Comm. Pure Appl. Math.�69�(2016),�no. 4,�734-744.

\bibitem{HaP83}  S. W. Hawking and D. N. Page; {\em  Thermodynamics of Black Holes In Anti-De
Sitter Space}, Comm. Math. Phys., \textbf{87} (1983) 577-588.

%\bibitem{Hi74} N. Hitchin, {\em Compact four-dimensional Eintein manifolds}, J. Diff. Geom. \textbf{9} (1974), 435-441.


%\bibitem{He96} E. Hebey; {em  Sobolev spaces on Riemannian manifolds},  Lecture Notes in Mathematics, Research Monograph, Springer-Verlag, Volume 1635,  1996.

\bibitem{HeK98}  M. Henningson and K. Skenderis; {\em The holographic Weyl anomaly}, J. High Ener. Phys., \textbf{07} (1998), 023, hep-th/9806087; Holography and the Weyl anomaly,
hep-th/9812032.

%\bibitem{HeK} M. Henningson and K. Skenderis, Weyl anomaly for Wilson surfaces, hep-th/9905163.

%\bibitem{Kr82} N.V. Krylov; {\em Boundedly nonhomogeneous elliptic and
%parabolic
%equations}, Izv. Akad. Nak. SSSR Ser. Mat., 46 (1982), 487-523;
%English transl. in Math. USSR Izv., 20 (1983), 459-492.
\bibitem{Ho} C. W. Hong; {\em A best constant and the Gaussian curvature}, Proc. Amer. Math. Soc., \textbf{97} (1986), 737-747.


\bibitem{LeB86} C. Lebrun: {\em On the topology of self-dual 4-manifolds}, Proc. Amer. Math. Soc., \textbf{98} (1986), 637-640.

\bibitem{LNN97} C. Lebrun, S. Nayatani, and T. Nitta; {\em Self-dual manifolds with positive Ricci curvature}, Math. Zeit. \textbf{224} (1997), 49-63.


\bibitem{Lee95} J. M. Lee; {\em The spectrum of an asymptotically hyperbolic Einstein manifold}, Comm. Anal. Geom., \textbf{3} (1995), no. 1-2, 253-271.

\bibitem{LiQS17} G. Li, J. Qing, and Y. Shi;  {\em Gap phenomena and curvature estimates for conformally compact Einstein manifolds}, 
Trans. Amer. Math. Soc., \textbf{369} (2017), 4385-4413.

\bibitem{Ma98} J. Maldacena; {\em The large N limit of superconformal field theories and super-
gravity};  Adv. Theor. Math. Phys. 2 (1998), 231-252, hep-th/9711200.

\bibitem{Mar98} C. Margerin; {\em A sharp characterization of the
smooth 4-sphere in curvature forms}, Comm Anal. Geom., \textbf{6} (1998), no. 1, 21-65.

\bibitem{Marq12} F.C. Marques; {\em Deforming three-manifolds with positive scalar curvature}, Ann. Math, \textbf{176} (2012), 815-863.


\bibitem{M71} J. Moser; {\em A Sharp form of an inequality by N.
Trudinger}, Indiana Math. J., \textbf{20} (1971), 1077-1091.

\bibitem{M71b} J. Moser; {\em On a non-linear problem in differential
geometry}, Dynamical Systems, (Proc. Sympos. Univ. Bahia, Salvador, 1971)
273-280. Academic Press, New York 1973.

\bibitem{On} E. Onofri; {\em On the positivity of the effective action
in a theory of random surfaces}, Comm. Math. Phys., \textbf{86} (1982), 321-326.

\bibitem{OPS88} B. Osgood, R. Phillips, and P. Sarnak; {\em Extremals
of determinants of Laplacians}, J. Funct. Anal.,  \textbf{80} (1988), 148-211.

\bibitem{OPS88b} B. Osgood, R. Phillips, and P. Sarnak; {\em Compact
 isospectral sets of surfaces}, J. Funct. Anal., \textbf{80} (1988), 212-234.

%\bibitem{Ok-1} K. Okikiolu: {\em The Campbell-Hausdorff theorem
%for elliptic operators and a related trace formula},
%Duke Math. J., 79 (1995) 687-722.

\bibitem{Ok01} K. Okikiolu; {\em Critical metrics for the determinant of
the Laplacian in odd dimensions}, Ann. Math., \textbf{153} (2001), no. 2,
471-531.

\bibitem{P78} D. Page; {\em A compact rotating gravitational instanton}, Phys. Lett. B, \textbf{79} (1978), 235-238.

\bibitem{Pa83} S. Paneitz; {\em A quartic conformally covariant
differential operator for arbitrary pseudo-Riemannian manifolds},
Preprint, 1983.

\bibitem{Po} A. Polyakov; {\em Quantum geometry of Bosonic strings},
Phys. Lett. B, \textbf{103} (1981), 207-210.

%\bibitem{Pog71} A.V. Pogorelev; {\em The Dirichlet problem for the multidimensional analogue of the Monge-Ampere equation}, Dokl. Acad. Nank. %SSSR, 201(1971), 790-793. In English translation, Soviet Math. Dokl. 12 (1971), 1227-1231.


\bibitem{Q03} J. Qing; {\em On the rigidity for conformally compact Einstein manifolds}, Int. Math. Res. Not., (2003), no. 21, 1141-1153.

\bibitem{S84} R. Schoen; {\em Conformal deformation of a Riemannian
metric to constant scalar curvature}, J. Diff. Geom., \textbf{20} (1984),
479-495.

\bibitem{S87} R. Schoen; {\em Variational Theory for the Total Scalar Curvature Functional
for Riemannian Metrics and Related Topics}, Springer Lecture notes 1365, (1987), 120-154.

\bibitem{T00} G. Tian; {\em  Canonical Metrics in K\"ahler Geometry},  Berlin: Birkh\"auser (2000).

%\bibitem{TV05a} G. Tian and J. Viaclovsky, \emph{Bach-flat asymptotically locally Euclidean metrics}, Invent. Math. \textbf{160} (2005), 357-415.

%\bibitem{TV05b} G. Tian and J. Viaclovsky, \emph{Moduli spaces of critical Riemannian metrics in dimension four}, Adv. Math. \textbf{196} (2005), 346-372.


%\bibitem{Tr} N. Trudinger; {\em On embedding into Orlicz spaces and some
% applications}, J. Math. Mech., 17 (1967), 473-483.

\bibitem{Tr68} N. Trudinger; {\em Remarks concerning the conformal deformation of Riemannian structure on compact manifolds}, Ann. Scuolo Norm. Sip. Pisa, \textbf{22} (1968), 265-274.


\bibitem{W98} E. Witten; {\em Anti-de Sitter space and holography};  Adv. Theor. Math. Phys., \textbf{2}
(1998), 253-290, hep-th/9802150.

%\bibitem{V00} J. Viaclovsky; {\em Conformal geometry, Contact geometry and the Calculus of Variations}, Duke Math. J., 101 (2000), no.2, 283-316.

\bibitem{Y60} H. Yamabe; {\em On a deformation of Riemannian structures on
compact manifolds}, Osaka Math. J., \textbf{12} (1960), 21-37.

%\bibitem{Y92} D. Yang;  {\em $L^p$ pinching and compactness theorems for compact Riemannian manifolds}, Forum Math. \textbf{4} (1992), 323-333.


\end{thebibliography}
\end{document}